# MINIMA IN BRANCHING RANDOM WALKS

By Louigi Addario-Berry[1] and Bruce Reed

*Université de Montréal, and McGill University and INRIA*

Given a branching random walk, let $M_n$ be the minimum position of any member of the $n$th generation. We calculate $\mathbf{E}M_n$ to within $O(1)$ and prove exponential tail bounds for $\mathbf{P}\{|M_n - \mathbf{E}M_n| > x\}$, under quite general conditions on the branching random walk. In particular, together with work by Bramson [*Z. Wahrsch. Verw. Gebiete* **45** (1978) 89–108], our results fully characterize the possible behavior of $\mathbf{E}M_n$ when the branching random walk has bounded branching and step size.

**1. Introduction.** The object of study in this paper is a supercritical branching random walk (or *tree-indexed random walk*)—which we view as a Galton–Watson tree $T$ with root $r$ and offspring distribution $B$, each of whose edges $e$ has been augmented (or *labeled*) with an independent copy $X_e$ of a random variable $X$ (which we call the *step size*). For the formal details of a probabilistic construction of branching random walks, see, for example, Harris [23]. When $X$ is nonnegative, this is often called an age-dependent branching process.

We assign to each node $v$ of $T$ the "displacement" (or *label*) $S_v$ which is the sum of the edge labels on the path from $r$ to $v$. The *depth* of a node $v \in T$ is the number of edges on the path from $r$ to $v$. One of the most well-studied parameters associated with branching random walks is the *minimum after $n$ steps*, that is, the minimum value of $S_v$ over all nodes $v$ having depth $n$ in $T$; we denote this quantity $M_n$. When $X$ is nonnegative, $M_n$ may be viewed as the time at which the first element of the $n$th generation is born. For several special choices of branching random walk, $M_n$ also turns out to be closely linked to the performance of data structures that arise in computer science; this connection is explored in depth by Devroye [18].

Received December 2007; revised July 2008.

[1]Supported in part by an NSERC discovery grant.

The starting point for our research is the Biggins–Hammersley–Kingman theorem, which provides a law of large numbers for the minimum displacement $M_n$ of a branching random walk after $n$ steps. (We set $M_n = \infty$ if the process does not survive for $n$ generations.) In a sequence of papers, each building on the results of the last, Hammersley [22], Kingman [25] and Biggins [6] showed that under suitable conditions on the exponential moment of $X$, there is a finite constant $\gamma$ such that, conditional on the survival of the branching process,

$$M_n/n \to \gamma \qquad \text{almost surely.}$$

When Hammersley initiated this research into the first-order behavior of $M_n$, he posed several questions to which complete answers remain unknown. In particular, he asked when more detailed information about $M_n - \gamma n$ than that given by the above law of large numbers can be found, about the expectation of $M_n$, and about whether the higher centralized moments of $M_n$ (and particularly the variance) are bounded. Our goal in this paper is to provide answers to Hammersley's questions for a broad range of branching random walks. Before we state our results, however, we briefly discuss some previous work on this and related problems.

In a remarkable paper, Bramson [7] derived extremely precise information about the maximum displacement of *branching Brownian motion*. In this model, an initial particle starts at position 0 on the real line and begins a standard Brownian motion with variance 1. The particle decays according to an exponential mean 1 clock; when the clock goes off, the particle splits into two, each of which continues an independent Brownian motion and each of which independently decays according to an exponential mean 1 clock. This process is continued forever. Bramson studied the maximum displacement $M_t^{\mathrm{br}}$ of any particle after time $t$ (of course, his results apply immediately to the minimum, by symmetry). He showed that the median $m_t^{\mathrm{br}}$ of $M_t^{\mathrm{br}}$ satisfies

$$(1) \qquad m_t^{\mathrm{br}} = \sqrt{2}t - \frac{3 \log t}{2} + c + o(1),$$

where $c$ is a fixed constant, and additionally showed that $M_t^{\mathrm{br}} - m_t^{\mathrm{br}}$ converges in distribution. (In this paper, by log we always mean the natural logarithm.) It is worth noting that Lifshits [26] has apparently (see Hu and Shi [24]) provided an example of a branching random walk for which $M_n$ minus its median is tight but does not converge in distribution, so in general a result as strong as Bramson's will not hold in the branching random walk setting.

Since Bramson's work, it has been expected that for at least some branching random walks, the median $m_n$ of $M_n$ should exhibit similar behavior. However, such results have been slow in coming. McDiarmid [27] proved (among other results) that for a wide class of branching random walks, $M_n - \gamma n = O(\log n)$ almost surely. Bachmann [3, 4] studied the *tightness* of the family $\{M_n - m_n : n \in \mathbb{N}\}$, where $m_n$ is the median of $M_n$, showing that the aforementioned family is tight when $X$ has logarithmically concave distribution function. (We have not stated all the conditions on $X$ and $B$ required for his result.) More recently, Bramson and Zeitouni [9, 10] also proved tightness of the family $\{M_n - m_n : n \in \mathbb{N}\}$, under a strong local condition on the lower tail of $X$.

In the course of studying the height of random binary search trees (RB-STs), Reed [28] proved that for a branching random walk with deterministic binary branching and exponential mean 1 step size, there exist explicit constants $\alpha > 0$, $\beta > 0$ such that

$$m_n = \alpha n + \beta \log n + O(1), \tag{2}$$

the first such result for any branching random walk; Reed additionally proved exponential tail bounds for $\mathbf{P}\{|M_n - m_n| > x\}$. We remark that the difference in sign of the logarithmic terms in (1) and (2) is because the former equation bounds a maximum, whereas the latter equation bounds a minimum. Drmota [21] has also published an alternative proof of Reed's result, using completely different techniques. Chauvin and Drmota [13] have since extended Drmota's techniques to show exponential tail bounds for $M_n$ around $m_n$ when (among other conditions) the step size $X$ satisfies what they call an *intersection property*. In particular, their result implies that such tail bounds hold when $X$ has logarithmically concave distribution function. In a series of papers, Broutin and Devroye [11], Broutin, Devroye and McLeish [12] have recently extended the Hammersley–Kingman–Biggins theorem to a wide range of branching random walks in which dependency is allowed among the steps and in which the "cutoff" at which the minimum is measured is not deterministically $n$, but is determined by a second family of random variables coexisting on the edges of the Galton–Watson tree $T$.

The results and techniques of Bramson [7] and Reed [28], in particular, serve as guidance and inspiration for the work of this paper. For our approach, the asymptotic behavior of $M_n$ and of its median $m_n$ are best characterized in terms of the behavior of the *logarithmic moment generating function* (LMGF) $\Lambda$. For a real random variable $X$, we define the LMGF $\Lambda = \Lambda_X$ of $X$ by

$$\Lambda(t) = \Lambda_X(t) = \log \mathbf{E}\{e^{tX}\}.$$

To better understand the utility of the function $\Lambda$ in studying $M_n$, we first recall Chernoff's bounding technique: if $S_n = \sum_{i=1}^n X_i$ is a sum of $n$ independent copies of $X$, then for any $c \in \mathbb{R}$ and $u > 0$, by using Markov's inequality we have

$$\mathbf{P}\{S_n \leq cn\} = \mathbf{P}\{e^{-uS_n} > e^{-ucn}\} \leq \frac{\mathbf{E}\{e^{-uS_n}\}}{e^{-ucn}} = [\mathbf{E}\{e^{-u(X-c)}\}]^n,$$

so

$$(3) \qquad \mathbf{P}\{S_n \leq cn\} \leq \left[\inf_{u>0} \mathbf{E}\{e^{u(X-c)}\}\right]^n.$$

If there is $u$ such that $\mathbf{E}\{e^{u(X-c)}\} < (\mathbf{E}B)^{-1}$ (where $B$ is the branching distribution) then it follows that $\mathbf{P}\{S_n \leq cn\} \leq (\mathbf{E}B)^{-n}$, and it is fairly easy to show, using this fact and straightforward bounds on the growth of $T$, that $\liminf_{n\to\infty} M_n/n > c$. (To show that in this case $\lim_{n\to\infty} M_n/n$ actually exists and to determine the limit takes more work.)

The optimal choice for $u$ in (3) is that for which $\Lambda'(u) = c$, as may be informally seen by differentiating $\mathbf{E}\{e^{u(X-c)}\}$ with respect to $u$. Choosing $u$ in this fashion, and writing $\Lambda'(u)$ in place of $c$, gives the bound

$$(4) \qquad \mathbf{P}\{S_n \leq \Lambda'(u)n\} \leq [\mathbf{E}\{e^{uX}\} \cdot e^{-u\Lambda'(u)}]^n = e^{-n(u\Lambda'(u)-\Lambda(u))}.$$

It turns out that the upper bound given by (4) is almost tight; this is the substance of the "exact asymptotics for large deviations" first proved by Bahadur and Rao [5], and is the reason that the behavior of $\Lambda$ is key to our investigation. As our proof leans strongly upon the exact asymptotics for large deviations, we now take the time to formally introduce this result.

Given $X$, $S_n$ and $\Lambda$ as above, let $D_\Lambda$ be the set of values $t$ for which $\Lambda(t)$ is finite, and let $D_\Lambda^o$ be the interior of $D_\Lambda$. We say $X$ is a *lattice random variable* with period $p > 0$ if there is a constant $x$ such that $(X-x)/p$ is almost surely integer-valued, and $p$ is the largest real number for which this holds. We now state the asymptotic estimates for large deviations as they appear in Dembo and Zeitouni [17] (as Theorem 3.7.4):

THEOREM 1 (Bahadur and Rao [5]).   *Let $S = \{S_n\}_{n\in\mathbb{N}}$ be a random walk with step size $X$, and define $\Lambda$ and $D_\Lambda^o$ as above. Choose any $t \in D_\Lambda^o$ with $t < 0$ and any function $g(n)$ tending to plus infinity with $n$; then if $X$ is nonlattice then for any $a \in R$ with $a \leq \sqrt{n}/g(n)$,*

$$(5) \qquad \mathbf{P}\{S_n \leq \Lambda'(t)n - a\} = (1 + o(1)) \cdot \frac{e^{at-n(t\Lambda'(t)-\Lambda(t))}}{\sqrt{\Lambda''(t) \cdot 2\pi n}},$$

*uniformly over $a$ in the above range. Furthermore, if $X$ is a lattice random variable with period $p$, then for any $a$ which is a multiple of $1/p$ and for which $a \leq \sqrt{n}/g(n)$, the same result holds.*

In fact, this theorem is stated with $a$ constant in [17], but the proof in that book yields without modification the above formulation. Theorem 1 contains two cases: the values of $a$ for which (5) holds depend on whether $X$ is lattice or nonlattice. Technically, this fact necessitates a case analysis in our proof, but as the two cases are virtually identical, we do not bother with both of them. *For the remainder of the paper, we assume that the step size $X$ is nonlattice.*

We will shortly have use of the following corollary of Theorem 1. Let $f(t) = t\Lambda'(t) - \Lambda(t)$.

COROLLARY 2. *For $0 \geq t \in D_\Lambda^0$, $f(t)$ is infinitely differentiable and strictly decreasing. Furthermore, if $X$ is bounded from below then*

$$\lim_{t \to -\infty} f(t) = \log(1/\mathbf{P}\{X = \operatorname{ess\,inf} X\})$$

*(where we interpret the right-hand side as $\infty$ if $\mathbf{P}\{X = \operatorname{ess\,inf} X\} = 0$).*

PROOF. The facts that $f(t)$ infinitely differentiable in $D_\Lambda^0$ and is strictly decreasing when $t \leq 0$ follow from Dembo and Zeitouni [17], Lemma 2.2.5 and Exercise 2.2.24. To prove the second assertion, first assume without loss of generality that $\operatorname{ess\,inf} X = 0$. Next note that $\Lambda'(t) = \mathbf{E}\{Xe^{tX}\}/\mathbf{E}\{e^{tX}\} > 0$ for all $t \neq 0$. In particular, for all $t < 0$, by Theorem 1 we thus have

$$\mathbf{P}\{X = 0\}^n \leq \mathbf{P}\{S_n \leq \Lambda'(t)n\} = (1 + o(1))\frac{e^{-nf(t)}}{\sqrt{\Lambda''(t) \cdot 2\pi n}},$$

so we must have $f(t) \leq \log(1/\mathbf{P}\{X = 0\})$. As $t < 0$ was arbitrary it follows that $\lim_{t \to -\infty} f(t) \leq \log(1/\mathbf{P}\{X = 0\})$.

To see that in fact equality holds, we first show that $\lim_{t \to -\infty} \Lambda'(t) = 0$. Note that for any fixed $\varepsilon > 0$ and $t < 0$, $\mathbf{E}e^{tX} \geq e^{t\varepsilon/3}\mathbf{P}\{X < \varepsilon/3\}$ and $\mathbf{E}\{e^{tX/2}\mathbf{1}_{[X \geq \varepsilon]}\} \leq e^{t\varepsilon/2}\mathbf{P}\{X \geq \varepsilon\}$, so

$$(6) \qquad \frac{\mathbf{E}\{e^{tX/2}\mathbf{1}_{[X \geq \varepsilon]}\}}{\mathbf{E}e^{tX}} \leq e^{t\varepsilon/6}\frac{\mathbf{P}\{X \geq \varepsilon\}}{\mathbf{P}\{X < \varepsilon/3\}} \to 0,$$

as $t \to -\infty$. Next, given $\varepsilon$ as above, for all sufficiently large negative $t$ and all $x \geq \varepsilon$ we have $xe^{tx} \leq \varepsilon e^{tx/2}$, so

$$\mathbf{E}\{Xe^{tX}\} = \mathbf{E}\{Xe^{tX}\mathbf{1}_{[X < \varepsilon]}\} + \mathbf{E}\{Xe^{tX}\mathbf{1}_{[X \geq \varepsilon]}\}$$

$$\leq \varepsilon(\mathbf{E}\{e^{tX}\mathbf{1}_{[X < \varepsilon]}\} + \mathbf{E}\{e^{tX/2}\mathbf{1}_{[X \geq \varepsilon]}\}).$$

It follows that for such $t$,

$$\frac{\mathbf{E}\{Xe^{tX}\}}{\mathbf{E}e^{tX}} \leq \varepsilon\left(1 + \frac{\mathbf{E}\{e^{tX/2}\mathbf{1}_{[X \geq \varepsilon]}\}}{\mathbf{E}e^{tX}}\right) \to \varepsilon$$

as $t \to -\infty$, by (6). Since $\varepsilon > 0$ was arbitrary it follows that

$$\lim_{t \to -\infty} \Lambda'(t) = \lim_{t \to -\infty} \frac{\mathbf{E}\{Xe^{tX}\}}{\mathbf{E}e^{tX}} = 0,$$

as claimed. Now fix $N \in \mathbb{N}_{>0}$, and small $\delta > 0$, let $\varepsilon = \delta/N$ and let $t_\varepsilon < t_\delta < 0$ be chosen so that $\Lambda'(t_\delta) = \delta$, $\Lambda'(t_\varepsilon) = \varepsilon$ [such choices of $t_\varepsilon$ and $t_\delta$ exist for $\delta$ small enough since $\Lambda'(t) > 0$ for all $t \neq 0$ and $\lim_{t \to -\infty} \Lambda'(t) = 0$]. Then for any $t \leq t_\varepsilon$, $f(t) \leq f(t_\varepsilon)$. Furthermore, if we are to have $S_n \leq \varepsilon n$ then at most $n/N$ of the $X_i$ can satisfy $X_i > \delta$. Letting $p_\delta = \mathbf{P}\{X \leq \delta\}$, it follows that

$$(1 + o(1)) \frac{e^{-nf(t_\varepsilon)}}{\sqrt{\Lambda''(t_\varepsilon) \cdot 2\pi n}} = \mathbf{P}\{S_n \leq \varepsilon n\} \leq \mathbf{P}\left\{\text{Bin}(n, p_\delta) \geq n\left(1 - \frac{1}{N}\right)\right\},$$

so

$$\lim_{t \to -\infty} f(t) \geq f(t_\varepsilon) \geq \frac{1}{n} \log(1/\mathbf{P}\{\text{Bin}(n, p_\delta) \geq n(1 - 1/N)\}).$$

Since $\delta > 0$ and $N \in \mathbb{N}$ were arbitrary and $\lim_{\delta \to 0} p_\delta = \mathbf{P}\{X = 0\}$, it then follows by standard binomial estimates that

$$\lim_{t \to -\infty} f(t) \geq \log(1/\mathbf{P}\{X = 0\}). \qquad \square$$

Using the above exact asymptotics, it turns out that one can prove an analogue of Bramson's result for a wide range of branching random walks. The primary result of this paper is:

THEOREM 3. *Consider a supercritical branching random walk with branching distribution $B$ and nonconstant step size $X$, and suppose that the following conditions hold:*

(I) *there exists an integer $d \geq 2$ such that $\mathbf{P}\{B \leq d\} = 1$;*
(II) *there exists a real number $u > 0$ such that $\mathbf{E}e^{uX} < \infty$; and*
(III) *there exists $t \in D_\Lambda^o$, $t < 0$, such that $t\Lambda'(t) - \Lambda(t) = \log(\mathbf{E}B)$.*

*Let $\mathcal{S}$ be the event that the branching random walk survives. Then*

$$(7) \qquad \mathbf{E}\{M_n | \mathcal{S}\} = \Lambda'(t)n - \frac{3}{2t}\log n + O(1),$$

*and there exist constants $C > 0$, $\delta > 0$ depending only on $X$, such that for all $x \in \mathbb{R}$,*

$$(8) \qquad \mathbf{P}\{|M_n - \mathbf{E}\{M_n | \mathcal{S}\}| \geq x | \mathcal{S}\} \leq Ce^{-\delta x}.$$

A brief discussion of the conditions appearing in Theorem 3 are in order. Conditions (II) and (III) control the positive and negative exponential moments of $X$, respectively. Condition (III) is perhaps the least intuitive. For example, if $X$ is bounded from below, then by Corollary 2, (III) simply requires that $\mathbf{P}\{X = \text{ess inf } X\} < 1/\mathbf{E}B$. Together with a result of Bramson [8], this shows that condition (III) can not be removed. In particular, Bramson [8] has shown that if $X$ is bounded from below and $\mathbf{P}\{X = \text{ess inf } X\} = \frac{1}{\mathbf{E}B}$, then $\mathbf{E}M_n = (\text{ess inf } X) \cdot n + O(\log \log n)$ (Dekking and Host [16] contains further information about behavior of $M_n$ in this case). The following theorem highlights another situation in which $\mathbf{E}M_n$ does not have the form $\alpha n + \beta \log n + O(1)$:

THEOREM 4. *Given a supercritical branching random walk satisfying conditions* (I) *and* (II), *above, if* $X$ *is bounded from below and additionally* $\mathbf{P}\{X = \operatorname{ess\,inf} X\} > 1/\mathbf{E}B$, *then*

$$\mathbf{E}\{M_n|\mathcal{S}\} = (\operatorname{ess\,inf} X) \cdot n + O(1)$$

*and there exist constants* $C' > 0$, $\delta' > 0$ *depending only on* $X$ *such that for all* $x > 0$, $\mathbf{P}\{M_n > (\operatorname{ess\,inf} X) \cdot n + x|\mathcal{S}\} \leq C'e^{-\delta' x}$.

Theorem 4 is not particularly difficult. It may be proved using the techniques developed by Hammersley, Kingman and Biggins in the course of proving the Hammersley–Kingman–Biggins theorem, and is also essentially contained within Dekking and Host [16]. We provide a proof that allows us to highlight, in a simplified setting, a technique we later use in proving Theorem 3. We remark that together with the aforementioned result of Bramson [8], Theorems 3 and 4 completely determine the possible behavior of $M_n$ in the case that $X$ is bounded from below.

If $X$ is not bounded from below, then condition (III) essentially imposes that the left tail of $X$ decays "sufficiently" (and in particular, at least exponentially) quickly. In particular, if $X$ is not bounded from below and all exponential moments $(X - \mathbf{E}X)_-$ are finite, then (III) will hold.

Condition (II) is essentially necessary for our result to hold; in general, a relaxation of condition (II) may require a corresponding relaxation of the upper tail bound in (8). In particular, it may be the case that

$$\mathbf{P}\{\text{the root } r \text{ has exactly one child}|\mathcal{S}\} > 0,$$

in which case it is not hard to see that the upper tail of $M_n$ can decay no more quickly than the upper tail of $X$.

Condition (I) is a requirement imposed by our particular use of the second moment method in the course of the proof. Independently of the current work, Hu and Shi [24] have shown that assuming bounded step size, assuming that (III) holds and that $\mathbf{E}\{B^{1+\varepsilon}\} < \infty$ for some $\varepsilon > 0$, then $(M_n - \Lambda'(n))/\log n \to (-3/2t)$ in probability on $\mathcal{S}$ but that, surprisingly, $(M_n - \Lambda'(n))/\log n$ does *not* converge almost surely. It seems very likely, particularly in view of the former result, that condition (I) should not be necessary for the results of Theorem 3 to hold. However, we would expect different behavior in the case $\mathbf{E}B = \infty$.

REMARK 1. Condition (III) allows us to apply Theorem 1 to obtain precise estimates for tail probabilities. However, there are certain situations in which such estimates are available without recourse to Theorem 1, for example, by direct computation. The key properties we require are the following: there exist constants $c \in \mathbb{R}$ and $a, a' > 0$ such that

$$\mathbf{P}\{S_n \leq cn - x\} = \Theta(e^{-ax} \cdot \mathbf{E}B^{-n} \cdot n^{-1/2}),$$

uniformly over all $x = o(\sqrt{n})$, and $\mathbf{P}\{S_n \leq cn - x\} = O(e^{-a'x} \cdot \mathbf{E}B^{-n})$ for all $x > 0$. Whenever $X$ satisfies these properties [and still assuming that (I) and (II) hold], we obtain an analogue of Theorem 3. More precisely, letting $x^* = x^*(n)$ solve $e^{-ax^*} \cdot n^{-1/2} = n$, Theorem 3 holds (and with an identical proof) if we replace (7) by the conclusion that $\mathbf{E}\{M_n | \mathcal{S}\} = cn + x^* + O(1)$. (The reason for this choice of $x^*$ will become clear over the course of the paper.) In particular, if $X = -\log(U)$ where $U$ is from the family of Beta distributions, such estimates are available by direct computation; this case is of particular interest due to its link with a variety of search trees arising in computer science (see, e.g., Chauvin and Drmota [13], Devroye [19], Drmota [21], Reed [28]).

1.1. *Our approach.* We observe that it clearly suffices to prove Theorem 3 when $\mathbf{E}X = 0$, and we shall hereafter assume this is the case. To begin to describe our approach, we observe that for $\lfloor (\mathbf{E}B)^n \rfloor$ *independent* random walks with step size $X$, the expected minimum value $\mathbf{E}M_n^{\mathrm{ind}}$ of any node $w$ at depth $n$ in this "forest of random walks" is within $O(1)$ of the smallest $m$ for which

$$\mathbf{E}\{|\{\text{nodes } w \text{ at depth } n \text{ such that } S_w \leq m\}|\} \geq 1.$$

In other words, $\mathbf{E}M_n^{\mathrm{ind}}$ is within $O(1)$ of $\inf\{m : \mathbf{P}\{S_w \leq m\} \geq 1/\lfloor (\mathbf{E}B)^n \rfloor\}$, a value we can easily derive to within $O(1)$ using Theorem 1. To get an idea why this is the case, choose $t \in D_\Lambda^o$, $t < 0$, and such that $t\Lambda'(t) - \Lambda(t) = \log \mathbf{E}B$. For any real number $a$ which is $o(\sqrt{n})$, by Theorem 1 and a union bound we have

$$\mathbf{P}\{M_n^{\mathrm{ind}} \leq \Lambda'(t)n - a\} \leq (1 + o(1)) \lfloor (\mathbf{E}B)^n \rfloor \cdot \frac{e^{at - n(t\Lambda'(t) - \Lambda(t))}}{\sqrt{\Lambda''(t) \cdot 2\pi n}}$$

$$= (1 + o(1)) \frac{\lfloor (\mathbf{E}B)^n \rfloor}{(\mathbf{E}B)^n} \frac{e^{at}}{\sqrt{\Lambda''(t) \cdot 2\pi n}}$$

$$= e^{at - (\log n)/2 - O(1)},$$

so $\mathbf{P}\{M_n^{\mathrm{ind}} \leq \Lambda'(t)n - a\}$ is exponentially small in $((\log n)/2t - a)$ when $a = o(\sqrt{n})$ and $a > 0$. Similarly, if $M_n^{\mathrm{ind}} > \Lambda'(t)n - a$ then each of the $\lfloor (\mathbf{E}B)^n \rfloor$ random walks $S$ must have $S_n \geq \Lambda'(t)n - a$. By the independence of the random walks and by Theorem 1, we have

$$\mathbf{P}\{M_n^{\mathrm{ind}} > \Lambda'(t)n - a\} = (1 - \mathbf{P}\{S_n \leq \Lambda'(t) - a\})^{\lfloor (\mathbf{E}B)^n \rfloor}$$

$$= \left(1 - (1 + o(1)) \frac{e^{at - (\log n)/2 - O(1)}}{(\mathbf{E}B)^n}\right)^{\lfloor (\mathbf{E}B)^n \rfloor}$$

$$\leq \exp\{-(1 + o(1))e^{at - (\log n)/2 - O(1)}\},$$

so $\mathbf{P}\{M_n^{\mathrm{ind}} \geq \Lambda'(t)n - a\}$ is doubly exponentially small in $(a - (\log n)/2t)$ when $a = o(\sqrt{n})$ and $a < 0$. These inequalities do not quite yield a bound on $\mathbf{E}M_n^{\mathrm{ind}}$, but they do show that an extremely high proportion of the probability mass of $M_n^{\mathrm{ind}}$ lies near $\Lambda'(t)n + (\log n)/2t$, and by combining these bounds with a Chernoff bound for the lower tail of $M_n^{\mathrm{ind}}$, it is not hard to show that in fact

$$\mathbf{E}\{M_n^{\mathrm{ind}}\} = \Lambda'(t)n - \frac{\log n}{2t} + O(1).$$

We denote by $m^* = m_n^*$ the quantity $\Lambda'(t)n - \log n/2t$, and call $m^*$ the *breakpoint*. Whenever we write $m^*$ without subscript we always mean $m_n^*$. By linearity of expectation, it follows easily from Theorem 1 that in the branching random walk,

$$\mathbf{E}\{|\{w \text{ at depth } n : S_w \leq m\}|\} \gg 1 \qquad \text{when } m - m^* \gg 1,$$

$$\mathbf{E}\{|\{w \text{ at depth } n : S_w \leq m\}|\} \ll 1 \qquad \text{when } m^* - m \gg 1.$$

Intuitively, then, the equation (7) can be understood by splitting the term $(3/2t)\log n$ into two pieces and writing

$$(9) \qquad \mathbf{E}\{M_n | \mathcal{S}\} = m^* - \frac{\log n}{t} + O(1);$$

the term $(\log n)/t$ must then be explained by the dependence in the branching random walk model. The bulk of the work of this paper is in understanding and explaining why this dependence should contribute a term of just this form.

At a high level, to explain this term, we shall apply the second moment method to bound the number of certain "special" nodes of $T$. We will introduce a notion of "goodness" of nodes; whether or not a node $v$ is good will depend only on the shape of the random walk from $r$ to $v$. We shall study the properties of the good subset $G$ of the nodes of $T$ at depth $n$, first showing that $|G|$ is tightly concentrated around its mean, and then, with some additional work, showing that $\mathbf{E}M_n$ is in fact within $O(1)$ of $\inf_m \mathbf{E}\{|\{w \in G : S_w \leq m\}|\} \geq 1$. To begin to make this more concrete, we first explain the key property that "good" nodes will satisfy.

1.2. *Leading nodes.* Given exchangeable random variables $X_1, \ldots, X_n$ with associated partial sums $S_1, \ldots, S_n$, we say that $S_n$ is *leading* (or that $S$ is leading after $n$ steps) if

$$(10) \qquad S_i \geq \mathbf{E}\{S_i | S_n\} \qquad \text{for all } i = 1, \ldots, n;$$

equivalently, if $S_i \geq S_n \cdot (i/n)$ for all $i = 1, \ldots, n$. (This terminology was introduced by McDiarmid [27].) Similarly, we say that $S_n$ is *strictly leading* if the inequality in (10) is strict for all $i = 1, \ldots, n$. Given a node $v$ at depth $n$ in $T$, we say that $v$ is a leading (resp. strictly leading) node if the random walk from the root to $v$ is leading (resp. strictly leading).

The intuition for why leading nodes are useful may be gleaned from examples (A) and (C) of the preceding section. If a node $v$ at depth $n$ is leading, then knowing that $S_v \leq m$ does not increase the *expected* number of nodes $w$ at depth $n$ with $S_w \leq m$ by too much. This fact allows us to use the *second moment method* (i.e., some variant of Chebyshev's inequality) to bound from below the probability that some such node exists when the expected number of such nodes is $\Omega(1)$. (This line of argument via Chebyshev's inequality is quite common in combinatorial settings; see, for example, Alon and Spencer [1], Chapter 4.)

It turns out that $\mathbf{E}\{M_n\}$ *is within $O(1)$ of the smallest value $m_n$ for which the expected number of leading vertices $v$ at depth $n$ with $S_v < m_n$ is at least* 1. Though it may not be immediately obvious, the assertion of the previous sentence is in fact equivalent to (7); we now explain this equivalence.

To bound the probability that a vertex is leading, we use a classic combinatorial technique introduced by Andersen [2] (and first used in the context of branching random walk by Devroye and Reed [20]), that we call a *rotation argument*, which involves considering cyclic permutations of the random variables $X_1, \ldots, X_n$; we will explain this approach in more detail in Section 3.

By this method, we will straightforwardly be able to show that for all $v \in N_n$ and for values $m \leq \mathbf{E}M_n$ that are not too far from the breakpoint $m^*$,

$$(11) \qquad \mathbf{P}\{S_v \leq m \text{ and } v \text{ is a leading node}\} = \Theta\left(\frac{\mathbf{P}\{S_v \leq m\}}{n}\right).$$

(We shall make this more formal in Section 3.) It follows that for such $m$,

$$\mathbf{E}\{|\{v \in N_n, S_v \leq m, v \text{ is leading}\}|\} = \Theta\left(\frac{\mathbf{E}\{|\{v \in N_n, S_v \leq m\}|\}}{n}\right)$$

$$(12) \qquad\qquad = \Theta\left(\frac{(\mathbf{E}B)^n}{n} \cdot \mathbf{P}\{S_v \leq m\}\right),$$

so we are in fact asserting that $\mathbf{E}M_n$ is within $O(1)$ of the smallest value $m$ for which, for nodes $v$ at depth $n$ in $T$, $\mathbf{P}\{S_v \leq m\} \geq (\mathbf{E}B)^n/n$. It follows immediately by Theorem 1 that $m$ is within $O(1)$ of $m^* + (\log n)/t$, in accordance with (7) and (9). We may then view the term $(\log n)/t$ as *the correction required in order to find a leading node.* One of the key steps in proving Theorem 3 will be to understand the shape of random walks that end in leading nodes, and in particular how much time such walks spend near their (conditional) means.

1.3. *A little notation and a basic fact.* Let $T_n$ be the subtree of $T$ consisting of all nodes of depth at most $n$. We may view $T_n$ as a subtree of a rooted, labeled $d$-ary tree $T_n'$ with $n$ levels and with root $r$, in the following manner. For a given node $v$ of $T_n'$ with children $w_1, \ldots, w_d$, let $B_v$ be a copy of the offspring random variable $B$, and let $X_{v,1}, \ldots, X_{v,B_v}$ be independent copies of $X$. Let $\sigma_v$ be a uniformly random permutation of $\{1, \ldots, d\}$; we assign label $X_{v,i}$ to edge $vw_{\sigma_v(i)}$ for $1 \le i \le B_v$ and assign label $\infty$ to edge $vw_{\sigma_v(i)}$ for $B_v < i \le d$. (This permutation gives the model a useful symmetry property that will be explained shortly.) We repeat this procedure for all nodes of $T_n'$, and label each node $v$ with the sum $S_v$ of all edge labels on the path from $r$ to $v$ (setting $S_v = \infty$ if any of these labels are infinite). With this labeling, the subtree of $T_n'$ induced by nodes $v$ with $S_v < \infty$ is distributed precisely as $T_n$, and we hereafter view $T_n = T_n(T_n')$ as a subtree of such a $d$-ary tree $T_n'$.

For $1 \le i \le n$ and for any $m \in \mathbb{R}$, let $N_i$ (resp. $N_i'$) be the set of nodes of $T_n$ (resp. $T_n'$) at depth $i$, and let $N_{i,m}$ be those nodes $v$ in $N_i$ with $S_v \le m$. To each $v \in T_n'$ we assign a label $S_v$ that is the sum of the edge labels $X_e$ on the path from $r$ to $v$—so it is possible that $S_v = \infty$ for some nodes $v$. We observe that for any integer $m \ge 0$, $M_n$ is the smallest $m \le +\infty$ for which $N_{n,m} \ne \varnothing$.

For the sake of our analysis, it will be useful to fix a distinguished path $P$ in $T_n'$ with nodes $r = v_0, v_1, \ldots, v_n$ (it may be useful to think of this path as running "along the left-hand side" of $T_n'$). Each node $v_i$ has one child $v_{i+1} = v_{i+1}^{(0)}$ in $P$—let its other children be called $v_{i+1}^{(1)}, \ldots, v_{i+1}^{(d-1)}$. Denote by $T^{i,j}$ the subtree of $T_n'$ rooted at $v_i^{(j)}$. Let $T_n^{i,j} = T^{i,j} \cap T_n$ (which may be empty), and let $N_{n,m}^{i,j}$ be the set of nodes of $N_{n,m}$ that are contained in $T_n^{i,j}$.

Given an automorphism $\alpha(T_n')$ of $T_n'$, we may view $\alpha$ as acting on the labels $X_{v,i}$ and on the permutations $\sigma_v$, by viewing the permutations $\sigma_v$ as fixed to the nodes of $T_n'$, the labels $X_{v,i}$ as fixed to the edges of $T_n'$, and both as being moved by the automorphism $a$. The presence of the permutations $\sigma_v$ ensures that for any automorphism $\alpha(T_n')$ of $T_n'$, the distribution of $T_n'$ (with both permutations and labels attached) is identical to that of $\alpha(T_n')$ (with both permutations and labels attached), and in particular, $\alpha$ induces a labeled isomorphism between $T_n(T_n')$ and $T_n(\alpha(T_n'))$. This symmetry, gained by the addition of the permutations, will greatly simplify later calculations. In particular, it yields the following fact.

FACT 5. *Let $E, F$ and $G$ be events in the $\sigma$-field generated by $T_n'$ with its labeling such that $E$ is invariant under automorphisms of $T_n'$ and there is some automorphism $\alpha$ of $T_n'$ such that $G = \alpha(F)$. Then $\mathbf{P}\{E|F\} = \mathbf{P}\{E|G\}$.*

For example, if $F$ is the event that $S_v \le m$ for some $m \in \mathbb{R}$ and some given $v \in N_n$, then $G$ could be the event that $S_{v'} \le m$ for any given $v' \in N_n$.

We will often use the phrase "by symmetry" in our arguments, rather than making explicit reference to Fact 5.

1.4. *Outline.* In Section 2, we prove Theorem 4, and in doing so introduce the idea of *amplification*, which also plays a role in the proof of Theorem 3. In Sections 3 and 4.1, we flesh out the high-level discussion of Reed's approach given above and explain some of the key ideas behind our proof of Theorem 3, particularly the upper bound on $M_n$ and the importance of leading nodes. In Section 4.2 we prove two key lemmas which allow us to control the shape of random walks conditioned on their value after time $n$. Finally, in Sections 4.3 and 5, respectively, we prove lower tail bounds and upper tail bounds on $M_n$ that together prove Theorem 3.

**2. Proof of Theorem 4.** To give an idea of how we will prove Theorem 4, we first consider the special case that the offspring distribution $B$ is deterministically $d$ and that $X$ is deterministically bounded from above, say $X \leq A$ for some constant $A$. In this case $\mathbf{P}\{S\} = 1$ so the conditioning in Theorem 4 vanishes. We also presume for simplicity that $\operatorname{ess\,inf} X = 0$. We consider the related branching process $T_0$ in which the set of children of a node is the set of its children in $T$ for which the displacement is 0. As $\mathbf{P}\{X = 0\} > 1/\mathbf{E}B = 1/d$, $T_0$ survives with positive probability, so there is a positive probability $p_0$ that $M_n = 0$ for every $n$.

Suppose that we want to bound the probability that $M_n$ is greater than $x$, for some given $x$. Since $X \leq A$ almost surely, if $n$ is at most $x/A$ then every node at depth $n$ has label at most $nA \leq x$, so $\mathbf{P}\{M_n > x\} = 0$. For larger $n$, we first observe that for any node $v$ at depth $\lfloor x/A \rfloor$, the tree rooted at $v$ whose nodes are the descendants $w$ of $v$ with $S_w = S_v$ is distributed precisely as $T_0$; we temporarily denote this tree $T_0(v)$.

The tree $T_0(v)$ survives with probability $p_0$, and if $T_0(v)$ survives then in particular, there is a descendent $w$ of $v$ at depth $n$ in $T(\mathcal{E})$ for which $S_w = S_v$, so $M_n \leq S_w = S_v \leq x$. It follows that

$$\mathbf{P}\{M_n > x\} \leq \mathbf{P}\left\{\bigcap_{v \text{ at depth } \lfloor x/A \rfloor} \{T_0(v) \text{ does not survive}\}\right\}.$$

Since the subtrees rooted at distinct nodes at depth $\lfloor x/A \rfloor$ are independent and there are $d^{\lfloor x/A \rfloor}$ such nodes, it follows that

$$\mathbf{P}\{M_n > x\} \leq \prod_{v \text{ at depth } \lfloor x/A \rfloor} \mathbf{P}\{T_0(v) \text{ does not survive}\}$$

$$(13) \qquad\qquad = (1 - p_0)^{d^{\lfloor x/A \rfloor}},$$

which in particular proves the tail bound of Theorem 4 and also immediately implies that $\mathbf{E}M_n = O(1)$. The key to the above line of reasoning is the idea of analyzing the subtrees of $T$ rooted at depth $\lfloor x/A \rfloor$ independently in order to strengthen our probability bound. We will hereafter refer to this technique as an *amplification* argument. McDiarmid [27] uses this idea in much the same fashion as above in his analysis of the minima of branching random walks; it also plays a key role in both Devroye and Reed [20] and Reed [28].

When $\mathbf{P}\{B = d\} < 1$ and $X$ is not necessarily bounded, we can not argue as straightforwardly as above. However, we still have that $\mathbf{P}\{T_0 \text{ survives}\} = p_0$ for some $p_0 > 0$. For $x \geq 0$ integers $d \geq 0$, we temporarily let $F_{x,d}$ be the set of nodes at depth $d$ with label at most $x$. Using an amplification argument just as we did in deriving (13) immediately yields that for any integers $c$ and $n$ with $c > 0$ and $n \geq d$,

$$\mathbf{P}\{M_n > x | |F_x| = c\} \leq (1 - p_0)^c. \tag{14}$$

So, to handle this case, we really need to control the distribution of the number of nodes at a given level of a supercritical branching process whose labels are not too large. To do so, we use the following result of McDiarmid ([27], Lemma 1), who showed that for any supercritical branching process there exist constants $\gamma_0 > 1$, $c_0 > 0$, and $\delta_0 > 0$ such that for all integers $i \geq 1$,

$$\mathbf{P}\{0 < |N_i| \leq \gamma_0^i\} \leq c_0 e^{-\delta_0 i}. \tag{15}$$

Fix $\gamma_0$, $c_0$ and $\delta_0$ as above, and define the function $\ell(x) = \lceil \log_{\gamma_0} x \rceil$. Since $\mathcal{S}$ occurs precisely if $|N_i| > 0$ for all $i$, for $n \geq \ell(x)$,

$$\mathbf{P}\{M_n > x, \mathcal{S}\} \leq \mathbf{P}\{0 < |N_{\lfloor \ell(x) \rfloor}| \leq x\} + \mathbf{P}\{M_n > x | |N_{\ell(x)}| > x\}$$
$$\leq c_0 e^{-\delta_0 x} + \sup_{x < k \leq d^{\ell(x)}} \mathbf{P}\{M_n > x | |N_{\ell(x)}| = k\}. \tag{16}$$

For fixed $k$ in the above range, we have

$$\mathbf{P}\{M_n > x | |N_{\ell(x)}| = k\}$$
$$= \sum_{i=0}^{k} (\mathbf{P}\{M_n > x | |N_{\ell(x)}| = k, |F_{x,\ell(x)}| = i\}$$
$$\qquad \cdot \mathbf{P}\{|F_{x,\ell(x)}| = i | |N_{\ell(x)}| = k\}) \tag{17}$$
$$\leq \sum_{i=0}^{k} (1 - p_0)^i \mathbf{P}\{|F_{x,\ell(x)}| = i | |N_{\ell(x)}| = k\}$$
$$\leq (1 - p_0)^{k/2} + \mathbf{P}\{|F_{x,\ell(x)}| \leq k/2 | |N_{\ell(x)}| = k\}.$$

Now fix $u > 0$ such that $\mathbf{E} e^{uX} = a < \infty$; such $u$ and $a$ exist by condition (II). We then have

$$\mathbf{P}\{v_{\ell(x)} \notin F_{x,\ell(x)} | v_{\ell(x)} \in N_{\ell(x)}\} = \mathbf{P}\{x < S_{v_{\ell(x)}} < \infty\}$$

$$\leq \frac{[\mathbf{E}e^{uX}]^{\ell(x)}}{e^{ux}}$$

(18)

$$= e^{(\log a)\lceil \log_{\gamma_0} x \rceil - ux}$$

$$\leq c_1 e^{-\delta_1 x},$$

for some $c_1 > 0$ and $\delta_1 > 0$ and for all $x > 0$. It follows that

$$\mathbf{E}\{|N_{\ell(x)} \setminus F_{x,\ell(x)}| \, | |N_{\ell(x)}| = k\} \leq c_1 e^{-\delta_1 x} \cdot k,$$

so by Markov's inequality

(19)     $$\mathbf{P}\{|F_{x,\ell(x)}| \leq k/2 \, | |N_{\ell(x)}| = k\} \leq 2c_1 e^{-\delta_1 x}.$$

Combining (16), (17) and (19), we thus have

(20)     $$\mathbf{P}\{M_n \geq x, \mathcal{S}\} \leq c_2 e^{-\delta_2 x},$$

for some $c_2 > 0$ and $\delta_2 > 0$ and all $x > 0$. Finally, since $T$ is supercritical, we have $\mathbf{P}\{\mathcal{S}\} = p > 0$. Letting $c_3 = c_2/p$, by (20) we thus have $\mathbf{P}\{M_n > x | \mathcal{S}\} \leq c_3 e^{-\delta_2 x}$ for $n \geq \ell(x)$. For $n < \ell(x)$, since given $\mathcal{S}$ there is at least one node in $N_n$, by symmetry and arguing as in (18) we have

$$\mathbf{P}\{M_n > x | \mathcal{S}\} \leq \mathbf{P}\{x < S_{v_n} < \infty\}$$
$$\leq c_1 e^{-\delta_1 x}.$$

This proves the exponential tail bounds of Theorem 4 and also shows that $\mathbf{E}\{M_n | \mathcal{S}\} = O(1) = \operatorname{ess\,inf} X + O(1)$. We remark that since we are assuming $\operatorname{ess\,inf} X = 0$, there is a shorter proof in the case that $n < \ell(x)$ since $M_{\ell(x)} \leq x$ implies $M_n \leq x$. We have given the above argument since it does not use the fact that $X$ is bounded from below, and we will appeal to it when proving Theorem 3.

**3. Typical leading nodes.** We recall that a node $v \in T_n$ is strictly leading if the random walk ending at $v$ stays strictly above its conditional expected value given $S_v$, that is, if it satisfies (10) with strict inequality. For $m \in \mathbb{R}$, we let

(21)     $$G_{n,m} = \{v \in N_n : m - 1 \leq S_v \leq m, S_v \text{ is strictly leading}\}.$$

We impose the requirement that $m - 1 \leq S_v$ as it gives us more precise control over the random walk; however, when $m$ is near $m^*$, Theorem 1 implies that $\mathbf{P}\{S_v \geq m - 1 | S_v \leq m\}$ is bounded away from zero. Furthermore, by replacing the constant $m - 1$ by $m - C$ for some large constant $C$, we could make this probability arbitrarily close to 1. Intuitively, therefore, we can think of $G_{n,m}$ as the set of "typical" leading nodes in $N_n$ with $S_v \leq m$. The following lemma is the promised formalization of (11).

LEMMA 6. *Given a random walk $S$ with steps distributed as $X$, and real numbers $a, c$ with $c > \mathbf{E}X \geq a$ and for which $\mathbf{P}\{\mathbf{E}X < X \leq c\} > 0$,*

$$\mathbf{P}\{S_n \leq an, S_n \text{ is strictly leading}\} \leq \frac{\mathbf{P}\{S_n \leq an\}}{n}$$

*and*

$$\mathbf{P}\{S_n \leq an, S_n \text{ is strictly leading}\} \geq \mathbf{P}\{\mathbf{E}X < X \leq c\} \cdot \frac{\mathbf{P}\{S_{n-1} \leq an - c\}}{n-1}.$$

PROOF. The proof of Lemma 6 is an adaptation of proofs from Andersen [2] and Reed [28]. For fixed $n$, and $n'$ with $n < n' \leq 2n$, let

$$S_{n'} = S_n + S_{n'-n} = \sum_{i=1}^{n'} X_{i(\bmod n)}.$$

We first note that if $S_n$ is strictly leading, then for any $j = 0, \ldots, n-1$, the random walk $S^{(j)}$ with $S_i^{(j)} = S_{j+i} - S_j$ for $i = 1, \ldots, n$ is not leading, since $S_n^{(j)} = S_n$ and

$$S_{n-j}^{(j)} = S_n - S_j < S_n\left(\frac{n-j}{n}\right).$$

More strongly, an identical argument shows that *at most* one of the random walks $S = S^{(0)}, S^{(1)}, \ldots, S^{(n)}$ is strictly leading. Furthermore, as the random variables $X_1, \ldots, X_n$ are independent and the event $\{S_n \leq an\}$ is fixed by permutations of $X_1, \ldots, X_n$, it follows that for all $j = 1, \ldots, n$,

$$\mathbf{P}\{S_n \leq an, S_n \text{ is strictly leading}\} = \mathbf{P}\{S_n \leq an, S_n^{(j)} \text{ is strictly leading}\},$$

and so

$$
\begin{aligned}
n \cdot \mathbf{P}&\{S_n \leq an, S_n \text{ is strictly leading}\} \\
&= \mathbf{P}\left\{\bigcup_{j=0}^{n-1}\{S_n \leq an, S_n^{(j)} \text{ is strictly leading}\}\right\} \\
&\leq \ \mathbf{P}\{S_n \leq an\},
\end{aligned}
$$

proving the upper bound of the lemma.

Next, for $i = 1, \ldots, n-1$, let $\hat{S}_i = S_{i+1} - X_1$. In order that $S_n \leq an$ and that $S_n$ is *strictly* leading, it suffices that:

- $\mathbf{E}X < X_1 \leq c$,
- $\hat{S}_{n-1} \leq an - c$ and
- $\hat{S}_{n-1}$ is leading.

As $X_1$ is independent from $\hat{S}_1, \ldots, \hat{S}_{n-1}$, it follows that

(22)
$$\mathbf{P}\{S_n \leq an, S_n \text{ is strictly leading}\}$$
$$\geq \mathbf{P}\{\mathbf{E}X < X_1 \leq c\} \cdot \mathbf{P}\{\hat{S}_{n-1} \leq an - c, \hat{S}_{n-1} \text{ is leading}\}.$$

For $j = 0, \ldots, n-2$, define the random walk $\hat{S}^{(j)}$ by $\hat{S}_i^{(j)} = \hat{S}_{j+i} - \hat{S}_j$, for $i = 1, \ldots, n-1$ (where, letting $n' = j + i$, if $n' > n-1$ then $\hat{S}'_n = \hat{S}_{n-1} + \hat{S}_{n'-(n-1)}$). Again we have that $\hat{S}_{n-1}^{(j)} = \hat{S}_{n-1}$ for all $j = 0, \ldots, n-2$. Furthermore, if $\hat{S}_{n-1}$ is not leading then, letting $j^*$ be an index for which $\hat{S}_{j^*} - \hat{S}_{n-1}(j^*/(n-1))$ is minimized, it follows immediately that for all $i = 1, \ldots, n-1$,

(23)
$$\hat{S}_i^{(j^*)} = \hat{S}_{i+j^*} - \hat{S}_{j^*}$$
$$= \hat{S}_{i+j^*} - \hat{S}_n\left(\frac{i+j^*}{n-1}\right) + \hat{S}_n\left(\frac{i+j^*}{n-1}\right) - \hat{S}_{j^*}$$
$$\geq \hat{S}_{j^*} - \hat{S}_n\left(\frac{j^*}{n-1}\right) + \hat{S}_n\left(\frac{i+j^*}{n-1}\right) - \hat{S}_{j^*}$$
$$= \hat{S}_n\left(\frac{i}{n-1}\right) = \hat{S}_n^{(j^*)}\left(\frac{i}{n-1}\right),$$

that is, $\hat{S}^{(j)}$ is leading at time $n-1$. Therefore, at least one of the random walks $\hat{S}^{(0)} = \hat{S}, \hat{S}^{(1)}, \ldots, \hat{S}^{(n-2)}$ is leading at time $n-1$. Since $X_2, \ldots, X_n$ are independent and the event $\{\hat{S}_{n-1} \leq an - c\}$ is fixed by permutations of $X_2, \ldots, X_n$, we thus have

(24)
$$\mathbf{P}\{\hat{S}_{n-1} \leq an - c\} = \mathbf{P}\left\{\hat{S}_{n-1} \leq an - c, \bigcup_{j=0}^{n-2}\{\hat{S}_{n-1}^{(j)} \text{ is leading}\}\right\}$$
$$\leq (n-1)\mathbf{P}\{\hat{S}_{n-1} \leq an - c, \hat{S}_{n-1} \text{ is leading}\}.$$

Combining (22) and (24), it follows that

$$\mathbf{P}\{S_n \leq an, S_n \text{ is strictly leading}\} \geq \mathbf{P}\{\mathbf{E}X < X_1 \leq c\} \cdot \frac{\mathbf{P}\{\hat{S}_{n-1} \leq an - c\}}{n-1},$$

proving the lower bound of the lemma. $\square$

By combining the above argument with the bounds from Theorem 1, the following lemma is immediate:

LEMMA 7. *Given any function $g(n)$ tending to infinity with $n$, for any $v \in N_n$ and any $m$ for which $|m^* - m| \leq \sqrt{n}/g(n)$,*

$$\mathbf{P}\{v \in G_{n,m}\} = \Theta\left(\frac{\mathbf{P}\{S_v \leq m\}}{n}\right),$$

*uniformly over all $m$ in the above range.*

We omit the proof of Lemma 7 as it is essentially identical to that of Lemma 6 (taking $a = m/n$ and choosing any fixed $c > \mathbf{E}X$ for which $\mathbf{P}\{\mathbf{E}X < X \le c\} > 0$; such $c$ exists since $X$ is nonconstant).

## 4. The upper bound.

4.1. *A warmup.* To prove an upper bound on $M_n$, we will eventually show that if $\mathbf{E}|G_{n,m}| = \Omega(1)$ then $\mathbf{P}\{|G_{n,m}| > 0\} = \Omega(1)$ [the definition of $G_{n,m}$ appears in (21)]. In order to demonstrate one of the key techniques of the lower bound in a simplified setting, we prove:

LEMMA 8. *For $m > m^*$ for which $m^* - m = o(\sqrt{n})$, if $\mathbf{E}|G_{n,m}| = \Omega(1/n^{5/2})$ then $\mathbf{P}\{|G_{n,m}| \ge 1\} = \Omega(1/n^{5/2}) = \Omega(1/m^{5/2})$.*

To prove Lemma 8, we use a version of the second moment method often called the Chung–Erdős inequality (see Chung and Erdős [15] and also Devroye and Reed [20] and Alon and Spencer [1], Chapter 2), which in our setting can be stated as follows: for any integer $i \ge 1$ and any random set $\mathcal{R} \subseteq N'_n$ (recall that $N'_n$ is the set of nodes at depth $n$ in $T'_n$),

$$(25) \qquad \mathbf{P}\{|\mathcal{R}| > 0\} \ge \frac{\mathbf{E}|\mathcal{R}|}{1 + \sup_{v \in N'_n} \mathbf{E}\{|\mathcal{R}| | v \in \mathcal{R}\}}.$$

We will apply (25) both here and later in the section. Recall that each $v_i$ on the distinguished path $v_0, \ldots, v_{n-1}$ has children $v_{i+1} = v_{i+1}^{(0)}, v_{i+1}^{(1)}, \ldots, v_{i+1}^{(d-1)}$ in $T'_n$. For $i = 0, \ldots, n$, $j = 0, \ldots, d-1$, and $m \in \mathbb{R}$ we hereafter denote

$$(26) \qquad N_{n,m}^{i,j} = \{v \in N_{n,m} : v \text{ is a descendent of } v_i^{(j)}\},$$

and define $G_{n,m}^{i,j}$ similarly.

PROOF OF LEMMA 8. By the symmetry of $T'_n$, $\mathbf{E}\{|G_{n,m}| | v \in G_{n,m}\}$ is identical for all $v \in N'_n$, so letting $\mathcal{R} = G_{n,m}$ in (25), we obtain

$$(27) \qquad \mathbf{P}\{|G_{n,m}| \ge 1\} \ge \frac{\mathbf{E}\{|G_{n,m}|\}}{1 + \mathbf{E}\{|G_{n,m}| | v_n \in G_{n,m}\}}.$$

By symmetry,

$$(28) \qquad \begin{aligned} \mathbf{E}\{|G_{n,m}|\} &= d^n \mathbf{P}\{v_n \in G_{n,m}\} \\ &= d^n \mathbf{P}\{v_n \in G_{n,m} | v_n \in N_n\} \mathbf{P}\{v_n \in N_n\}. \end{aligned}$$

Given that $v_n \in N_n$, $S_{v_n}$ is just a sum of $n$ i.i.d. random variables distributed as $X$, so by Lemma 7 and Theorem 1,

$$
\begin{aligned}
\mathbf{P}\{v_n \in G_{n,m} | v_n \in N_n\} &= \Theta\left(\frac{\mathbf{P}\{S_{v_n} \leq m | v_n \in N_n\}}{n}\right) \\
&= \Theta\left(\frac{e^{-t(m-m^*)}}{n \cdot [\mathbf{E}B]^d}\right),
\end{aligned}
\tag{29}
$$

where $t < 0$ has been chosen such that $t\Lambda'(t) - \Lambda(t) = \log(\mathbf{E}B)$. Now, for any edge $e = vw$ of the tree $T_n'$,

$$
\begin{aligned}
\mathbf{P}\{X_e < \infty\} &= \sum_{i=0}^{d} \mathbf{P}\{X_e < \infty | B_v = i\}\mathbf{P}\{B_v = i\} \\
&= \frac{1}{d}\sum_{i=1}^{d} i\mathbf{P}\{B_v = i\} = \frac{\mathbf{E}B}{d}.
\end{aligned}
\tag{30}
$$

As the variables $\{B_v\}_{v \in T_n'}$ are independent, it follows that for distinct edges $e = vw$ and $f = xy$, $X_e$ and $X_f$ are independent unless $v = x$, that is, unless $w$ and $y$ are siblings. In particular, it follows from this independence and from (30) that

$$
\mathbf{P}\{v_n \in N_n\} = \mathbf{P}\left\{\bigcap_{i=1}^{n}\{X_{v_{i-1}v_i} < \infty\}\right\} = \left(\frac{\mathbf{E}B}{d}\right)^n.
\tag{31}
$$

By combining (28), (29) and (31), we obtain the bound

$$
\mathbf{E}|G_{n,m}| = \Theta\left(\frac{e^{-t(m-m^*)}}{n}\right).
\tag{32}
$$

[A brief remark: we will usually omit the sorts of arguments leading to (31) when such calculations arise in later proofs; we have included them once for completeness.] We now turn to bounding $\mathbf{E}\{|G_{n,m}||v_n \in G_{n,m}\}$. By symmetry, we have

$$
\begin{aligned}
\mathbf{E}\{|G_{n,m}||v_n \in G_{n,m}\} &\leq \mathbf{E}\{|N_{n,m}||v_n \in G_{n,m}\} \\
&= 1 + \sum_{i=0}^{n-1}\sum_{j=1}^{d-1}\mathbf{E}\{|N_{n,m}^{i+1,j}||v_n \in G_{n,m}\} \\
&= 1 + (d-1)\sum_{i=0}^{n-1}\mathbf{E}\{|N_{n,m}^{i+1,1}||v_n \in G_{n,m}\}.
\end{aligned}
\tag{33}
$$

Since $v_n \in G_{n,m}$, for $i$ with $0 \leq i \leq n$ we have $S_{v_i} \geq (m-1)(i/n) \geq mi/n - 1$. It follows that for such $i$, $N_{n,m}^{i+1,1}$ is at most the number of descendants $v$ of $v_{i+1}^{(1)}$ for which $S_v - S_{v_i}$ is at most $m(n-i)/n + 1$. We recall that $m_i^* = \Lambda'(t)i - \log i/2t$. Let $\Delta_i = mi/n + 1 - m_i^*$ for $i = 1, \ldots, n$, and fix some distinguished descendant $v$ of $v_{i+1}^{(1)}$ at depth $n$. By the previous observation, linearity of expectation and symmetry we have

$$\mathbf{E}\{|N_{n,m}^{i+1,1}| | v_n \in G_{n,m}\}$$

$$\leq |N_n^{i+1,1}| \mathbf{P}\{S_v - S_{v_i} \leq m_{n-i}^* + \Delta_{(n-i)} | v_n \in G_{n,m}\}$$

$$(34) \qquad = d^{n-(i+1)} \cdot \left(\frac{\mathbf{E}B}{d}\right)^{n-(i+1)}$$

$$\cdot \mathbf{P}\{S_v - S_{v_i} \leq m_{n-i}^* + \Delta_{(n-i)} | v \in N_n^{i+1,1}\}$$

$$= (\mathbf{E}B)^{n-(i+1)} \cdot \Theta\left(\frac{e^{-t\Delta_{n-i}}}{(\mathbf{E}B)^{n-i}}\right) = \Theta(e^{-t\Delta_{n-i}}).$$

We now claim that for all $i = 1, \ldots, n$,

$$(35) \qquad e^{-t\Delta_{n-i}} = O(n^{1/2} e^{-t(m-m^*)}).$$

From (33)–(35), it follows that

$$(36) \qquad \mathbf{E}\{|G_{n,m}| | v_n \in G_{n,m}\} \leq 1 + O(n^{3/2} e^{-t(m-m^*)}).$$

Combining (32) and (36) we see that if $\mathbf{E}|G_{n,m}| = \Omega(1/n^{5/2})$ then

$$\mathbf{E}\{|G_{n,m}| | v_n \in G_{n,m}\} \leq 1 + O(n^{5/2} \cdot \mathbf{E}|G_{n,m}|),$$

which together with (28) proves the fact. To see that (35) holds, we write

$$\Delta_{n-i} = \frac{m(n-i)}{n} - m_{n-i}^* + 1$$

$$= \frac{m(n-i)}{n} - \Lambda'(t)(n-i) - \frac{\log(n-i)}{2t} + 1$$

$$= m - \frac{mi}{n} + \Lambda'(t)i - \Lambda'(t)n - \frac{\log n}{2t} + \frac{\log n - \log(n-i)}{2t} + 1$$

$$= m - m^* + \Lambda'(t)i - \frac{mi}{n} + \frac{\log n - \log(n-i)}{2t} + 1$$

$$< m - m^* + \frac{i}{n}(\Lambda'(t)n - m) + 1$$

$$\leq m - m^* - \frac{\log n}{2t} + 1,$$

the final inequality holding as $m \geq m^* = \Lambda'(t)n + \log n/(2t)$ (recall that $t$ is negative). Thus $\Delta_{n-i} \leq m - m^* - \log n/2t$, from which (35) follows immediately. $\square$

By Lemma 8 and a standard amplification argument, we could show that $\mathbf{P}\{|N_{n,m+a\log n}| \geq 1\} = \Omega(1)$ for some $a = O(1)$ and thereby show that $\mathbf{E}\{M_n|\mathcal{S}\} \leq m^* + a\log n$; we will not do so as we are headed toward stronger upper bounds, the key step of which is to strengthen Lemma 8. We now state this aim formally.

LEMMA 9.    *For all values of $m$ for which $m^* - m = o(\sqrt{n})$,*

(37)                $\mathbf{E}\{|G_{n,m}| | v_n \in G_{n,m}\} = \mathbf{E}|G_{n,m}| + O(1).$

Assume for the moment that Lemma 9 holds, and let

$$m' = m^* - \frac{\log n}{2t} = \Lambda'(t)n - \frac{3\log n}{2t}.$$

By Lemma 7 and linearity of expectation, we see that $m'$ is within $O(1)$ of the smallest value $m$ for which $\mathbf{E}|G_{n,m}| \geq 1$. Applying (25) with $\mathcal{R} = G_{n,m'}$, it follows immediately that $\mathbf{P}\{|G_{n,m'}| > 0\} = \Omega(1)$. Since $G_{n,m} \subseteq N_{n,m}$ for all $m$, it follows that $\mathbf{P}\{|N_{n,m'}| > 0\} = \Omega(1)$, so $\mathbf{P}\{M_n \leq m'\} = \Omega(1)$. By applying an amplification argument exactly as we did when proving Theorem 4, we immediately obtain exponential tail bounds for the upper tail of $M_n$:

COROLLARY 10.    *There exist constants $C_1 > 0, \delta_1 > 0$ such that for all $x > 0$,*

$$\mathbf{P}\{M_n \geq m' + x|\mathcal{S}\} \leq C_1 e^{-\delta_1 x}.$$

To obtain the upper tail bounds of Theorem 3, it thus remains to prove Lemma 9.

4.2. *The shape of random walks.*    The proof of Lemma 9 is based on establishing finer control over the behavior of a given random walk in $T$ than that given by Lemma 7. In this section, we prove the two lemmas that accomplish this. In Lemma 11, we bound the probability that a walk $S_0, S_1, \ldots, S_n$ is ever very far below its conditional mean given $S_n$; in Lemma 12, we bound the probability that such a random walk is ever *close to* its conditional mean, given that it is *never* very far below its conditional mean. These results are rather straightforward; the only technicalities are due to the fact that we must apply Theorem 1 in the course of the proofs. We now proceed to the details.

For $a \leq 0$, we say that $S$ *stays above $a$* up to time $n$, and write $\{S_n \ abo \ a\}$, if

$$S_i > S_n \cdot \frac{i}{n} + a \qquad \text{for all } i = 1, \ldots, n.$$

In particular, $\{S_n \ abo \ 0\}$ is simply the event that $S_n$ is strictly leading.

LEMMA 11. *Given any function $g(n)$ tending to plus infinity with $n$, any $m$ for which $|\Lambda'(t)n - m| = \sqrt{n}/g(n)$, and any $a \le -1$,*

$$\mathbf{P}\{S_n \ abo \ a | m-1 \le S_n \le m\} = O\left(\frac{|a|^9}{n}\right), \tag{38}$$

*uniformly over all $m$ and $a$ in the above ranges.*

PROOF. For simplicity we assume that $a$ is an integer, that $m = \Lambda'(t)n$, and that $g(n) \le \log n$; this eases the notational burden without changing the essence of the proof. The probability in (38) increases as $a$ decreases; it thus suffices to prove (38) when $|a|$ is at least some large fixed constant $C$. We assume $|a|$ is large enough that $|a| \le \sqrt{|a|^3}/\log(|a|^3) \le \sqrt{|a|}/g(|a|)$.

We remark that if $|a| > n^{1/9}$ then (38) holds trivially; we thus assume that $|a| \le n^{1/9}$. Let $n_1 = n + 2|a|^3$, and let $m_1 = \Lambda'(t)n_1$; as $|a| \le n^{1/9}$, we have $n_1 = \Theta(n)$. Let $S'$ be the random walk with $S'_i = S_{|a|^3+i} - S_{|a|^3}$ (the original walk "started at time $|a|^3$"). The walk $S'$ is distributed as $S$; we will show that

$$\mathbf{P}\{S'_n \ abo \ a | m-1 \le S'_n \le m\} = O\left(\frac{|a|^9}{n}\right), \tag{39}$$

which proves the lemma. We proceed by comparing the following two probabilities:

$$\mathbf{P}\{m_1 - 3 \le S_{n_1} \le m_1, S_{n_1} \ abo \ 0\} \quad \text{and} \quad \mathbf{P}\{m-1 \le S'_n \le m, S'_n \ abo \ a\}.$$

By Lemma 7,

$$\mathbf{P}\{m_1 - 3 \le S_{n_1} \le m_1, S_{n_1} \ abo \ 0\} = \Theta\left(\frac{\mathbf{P}\{m_1 - 3 \le S_{n_1} \le m_1\}}{n_1}\right)$$
$$= \Theta\left(\frac{\mathbf{P}\{m_1 - 3 \le S_{n_1} \le m_1\}}{n}\right). \tag{40}$$

By Theorem 1,

$$\mathbf{P}\{m_1 - 3 \le S_{n_1} \le m_1\} = \Theta\left(\frac{e^{-\Lambda'(t)n_1}}{\sqrt{n_1}}\right) = \Theta\left(\frac{e^{-2|a|^3\Lambda'(t)}e^{-\Lambda'(t)n}}{\sqrt{n}}\right)$$
$$= \Theta(e^{-2|a|^3\Lambda'(t)}\mathbf{P}\{m-1 \le S'_n \le m\}), \tag{41}$$

so, letting $E^*$ be the event $\{m_1 - 3 \le S_{n_1} \le m_1, S_{n_1} \ abo \ 0\}$, (40) and (41) yield that

$$\mathbf{P}\{E^*\} = \Theta\left(\frac{e^{-2|a|^3\Lambda'(t)}\mathbf{P}\{m-1 \le S'_n \le m\}}{n}\right). \tag{42}$$

We will show that

$$\mathbf{P}\{E^*\} = \Omega\left(\frac{e^{-2|a|^3 \Lambda'(t)}\mathbf{P}\{m-1 \leq S'_n \leq m, S'_n \; abo \; a\}}{|a|^9}\right). \tag{43}$$

Equation (39) follows immediately from (42) and (43); it thus remains to prove (43). We claim that for $E^*$ to occur it suffices that the following events occur:

(1) $\Lambda'(t)|a|^3 + |a| + 1 \leq S_{|a|^3} \leq \Lambda'(t)|a|^3 + |a| + 2$, which we denote $E_1$;
(2) $\Lambda'(t)|a|^3 - (|a|+3) \leq S_{n_1} - S_{n_1-|a|^3} \leq \Lambda'(t)|a|^3 - (|a|+2)$, which we denote $E_2$;
(3) $m - 1 \leq S'_n \leq m$, which we denote $E_3$;
(4) $S_{|a|^3} \; abo \; 0$;
(5) letting $S_j^{(1)} = S_{n_1-|a|^3+j} - S_{n_1-|a|^3}$ for $j = 0, 1, \ldots, |a|^3$, $S_{|a|^3}^{(1)} \; abo \; 0$; and finally,
(6) $S'_n \; abo \; a$.

Events (1)–(6) are depicted in Figure 1. Informally, (1) and (4) ensure that the first $|a|^3$ steps of the walk do not prevent $E^*$ from occurring, (2) and (5) do the same for the last $|a|^3$ steps of the walk, and (3) and (6) do likewise for the middle $n_1 - 2|a|^3 = n$ steps of the walk. [This is not quite the whole story; the "extra height" (between $|a|+1$ and $|a|+2$) gained in (1), as well as the extra height gained in (2), are required so that (3) and (6) can do their job.] By the independence of disjoint sections of the random walk, we thus have

$$\mathbf{P}\{E^*\} \geq \mathbf{P}\{E_1, E_2, E_3, S_{|a|^3} \; abo \; 0, S_{|a|^3}^{(1)} \; abo \; 0, S'_n \; abo \; a\}$$

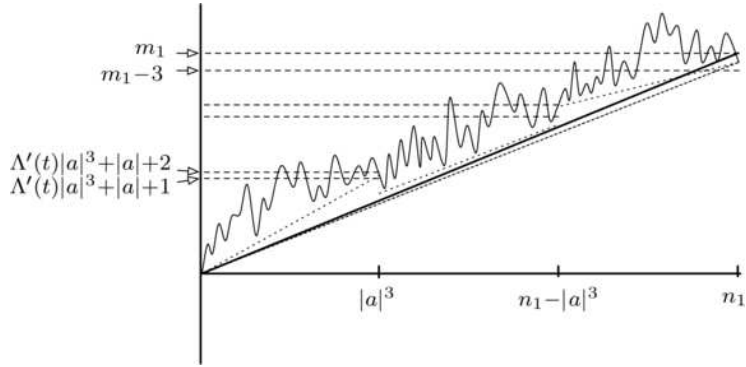

Fig. 1. *The events (1)–(6). For example, (1) occurs because after $|a|^3$ steps the walk is between the lowest two dashed horizontal lines, and (4) occurs because in the first $|a|^3$ steps the walk stays above the dashed line segment connecting its endpoints. The event $E^*$ occurs because after $n_1$ steps the walk is between $m_1 - 3$ and $m_1$ and has stayed above dashed line connecting its endpoints.*

(44)
$$= \mathbf{P}\{E_1, S_{|a|^3} \; abo \; 0\} \cdot \mathbf{P}\{E_2, S^{(1)}_{|a|^3} \; abo \; 0\} \cdot \mathbf{P}\{E_3, S'_n \; abo \; a\}.$$

As $|a| \le \sqrt{|a|^3}/g(|a|)$, we may apply Theorem 1 and Lemma 7 to bound $\mathbf{P}\{E_1, S_{|a|^3} \; abo \; 0\}$:

$$
\begin{aligned}
(45) \qquad \mathbf{P}\{E_1, S_{|a|^3} \; abo \; 0\} &= \Theta\left(\frac{\mathbf{P}\{E_1\}}{|a|^3}\right) \\
&= \Theta\left(\frac{e^{-|a|^3 \Lambda'(t)} e^{-t|a|}}{|a|^3 \sqrt{|a|^3}}\right) \\
&= \Theta\left(\frac{e^{-|a|^3 \Lambda'(t)} e^{-t|a|}}{|a|^{9/2}}\right).
\end{aligned}
$$

Likewise, we have

$$(46) \qquad \mathbf{P}\{E_2, S^{(1)}_{|a|^3} \; abo \; 0\} = \Theta\left(\frac{\mathbf{P}\{E_2\}}{|a|^3}\right) = \Theta\left(\frac{e^{-|a|^3 \Lambda'(t)} e^{t|a|}}{|a|^{9/2}}\right).$$

Combining (44)–(46) yields that

$$\mathbf{P}\{E^*\} = \Omega\left(\frac{e^{-2|a|^3 \Lambda'(t)}}{|a|^9} \mathbf{P}\{E_3, S'_n \; abo \; a\}\right),$$

which is precisely the claim in (43). □

Using Lemma 11, we can bound the conditional probability that $S$ spends much time *near* its mean, given that it is never much below its conditional mean up to time $n$. More precisely: for $0 \le k \le n$, let $b(n,k) = \min\{k, n-k\}$, and let $C_k$ be the event that $S_k \le S_n(k/n) + b(n,k)^{1/57}$ (that $S_k$ is "close" to its conditional mean; the meaning of "close" depends on how near $k$ is to one of the ends of the random walk). Given an integer $a \le -1$, let $B_a = \bigcup_{k=|a|^{57}}^{n-|a|^{57}} C_k$ (if $|a|^{57} > n/2$ then $B_a$ is the empty event). Then:

LEMMA 12. *Given any function $g(n)$ tending to plus infinity with $n$, for any $m$ for which $|\Lambda'(t)n - m| \le \sqrt{n}/g(n)$ and any integer $a \le -1$,*

$$(47) \qquad \mathbf{P}\{S_n \; abo \; a, B_a | m-1 \le S_n \le m\} = O\left(\frac{1}{n|a|^9}\right).$$

PROOF. As in the proof of Lemma 11, for the sake of readability we assume that $m = \Lambda'(t)n$ and that $g(k) \le \log k$ for all $k$. We may presume that $|a|^{57} < n/2$ or else the claim holds trivially. For any fixed $K$, for all $a$ with $|a| \le K$, the claim follows immediately from Lemma 11, so we may assume $|a|$ is larger than any fixed $K$; we henceforth assume $|a|$ is large enough that $|a| \leq \sqrt{|a|^{57}}/g(|a|)$.

Choose $k$ with $b(n,k) \geq |a|^{57}$. If $\{m-1 \leq S_n \leq m\}$ and $C_k$ are to both occur then necessarily

$$(48) \qquad S_k \leq \frac{mk}{n} + b(n,k)^{1/57} = \Lambda'(t)k + b(n,k)^{1/57}.$$

Similarly, if $\{m-1 \leq S_n \leq m\}$ and $\{S_n \; abo \; a\}$ are to both occur then necessarily

$$(49) \qquad S_k \geq \Lambda'(t)k + a - 1$$

(recall that $a$ is negative). More strongly, for $C_k$, $\{m-1 \leq S_n \leq m\}$ and $\{S_n \; abo \; a\}$ to *all* occur, there must be some integer $i$ with $a-1 \leq i \leq b(n,k)^{1/57}$ for which the following four events occur:

- $\Lambda'(t)k + i \leq S_k \leq \Lambda'(t)k + i + 1$, which we denote $A_{i,k}$;
- $\{S_k \; abo \; a - (i+1)\}$, which we denote $D_{i,k}$;

and letting $S'_j = S_{k+j} - S_j$ for $0 \leq j \leq n-k$,

- $\Lambda'(t)(n-k) - (i+1) \leq S'_{n-k} \leq \Lambda'(t)(n-k) - i$, which we denote $E_{i,k}$;
- $\{S'_{n-k} \; abo \; a-(i+1)\}$, which we denote $F_{i,k}$.

Restating the above using the names of these events, we have

$$C_k \cap \{S_n \; abo \; a\} \cap \{m-1 \leq S_n \leq m\} \subset \bigcup_{i=a-1}^{\lfloor b(n,k)^{1/57} \rfloor} A_{i,k} \cap D_{i,k} \cap E_{i,k} \cap F_{i,k}.$$
$$(50)$$

Since $B_a = \bigcup_{k=|a|^{57}}^{n-|a|^{57}} C_k$, It follows by a union bound that

$$(51) \qquad \begin{aligned} & \mathbf{P}\{B_a, S_n \; abo \; a, m-1 \leq S_n \leq n\} \\ & \leq \sum_{k=|a|^{57}}^{n-|a|^{57}} \sum_{i=a-1}^{\lfloor b(n,k)^{1/57} \rfloor} \mathbf{P}\{A_{i,k}, D_{i,k}, E_{i,k}, F_{i,k}\}. \end{aligned}$$

We will show that for each $k$ with $b(n,k) \geq |a|^{57}$, for all $i$ in the above range,

$$(52) \quad \mathbf{P}\{A_{i,i}, D_{i,k}, E_{i,k}, F_{i,k}\} = O\left(\mathbf{P}\{m-1 \leq S_n \leq m\} \cdot \frac{b(n,k)^{18/57} n^{1/2}}{(k(n-k))^{3/2}}\right).$$

We remark that the bound in (52) does not depend on $i$. Presuming for a moment that (52) holds, by dividing through by $\mathbf{P}\{m-1 \leq S_n \leq m\}$ in (51) and (52), we obtain

$$\mathbf{P}\{B_a, S_n \; abo \; a | m-1 \leq S_n \leq n\}$$

$$(53) \qquad = \sum_{k=|a|^{57}}^{n-|a|^{57}} \sum_{i=a-1}^{\lfloor b(n,k)^{1/57} \rfloor} \frac{\mathbf{P}\{A_{i,k}, D_{i,k}, E_{i,k}, F_{i,k}\}}{\mathbf{P}\{m-1 \le S_n \le m\}}$$

$$= O\left(\sum_{k=|a|^{57}}^{n-|a|^{57}} \sum_{i=a-1}^{\lfloor b(n,k)^{1/57} \rfloor} \frac{b(n,k)^{18/57} n^{1/2}}{(k(n-k))^{3/2}}\right).$$

Since $b(n,k)^{1/57} \ge |a|$, the inner sum in (53) has at most $2b(n,k)^{1/57} + 2$ identical terms. We thus have

$$\mathbf{P}\{B_a, S_n \ abo \ a | m-1 \le S_n \le n\} = O\left(\sum_{k=|a|^{57}}^{n-|a|^{57}} \frac{b(n,k)^{19/57} n^{1/2}}{(k(n-k))^{3/2}}\right)$$

$$= O\left(\sum_{k=|a|^{57}}^{\lceil n/2 \rceil} \frac{b(n,k)^{19/57} n^{1/2}}{(k(n-k))^{3/2}}\right)$$

$$= O\left(\frac{1}{n} \sum_{k=|a|^{57}}^{\lceil n/2 \rceil} \frac{k^{19/57}}{k^{3/2}}\right)$$

$$= O\left(\frac{1}{n} \sum_{k=|a|^{57}}^{\lceil n/2 \rceil} \frac{1}{k^{7/6}}\right)$$

$$= O\left(\frac{1}{|a|^9 \cdot n}\right),$$

which proves the lemma. It therefore remains to prove (52). To do so, we first write

$$\mathbf{P}\{A_{i,k}, D_{i,k}\} = \mathbf{P}\{A_{i,k}\} \cdot \mathbf{P}\{D_{i,k} | A_{i,k}\}.$$

Since $|i| \le \lfloor b(n,k) \rfloor^{1/57} \le k^{1/57} \le \sqrt{k}/g(k)$, and likewise $|a| \le b(n,k)^{1/57}$, we have $|a - (i+1)| \le 2\sqrt{k}/g(k) + 1$. By Lemma 11 (applied with a slightly different function $g$) it then follows that

$$(54) \qquad \mathbf{P}\{D_{i,k} | A_{i,k}\} = O\left(\frac{(|a| + |i| + 1)^9}{k}\right) = O\left(\frac{b(n,k)^{9/57}}{k}\right),$$

so

$$(55) \qquad \mathbf{P}\{A_{i,k}, D_{i,k}\} = O\left(\frac{\mathbf{P}\{A_{i,k}\} \cdot b(n,k)^{9/57}}{k}\right),$$

and an identical derivation yields

$$(56) \qquad \mathbf{P}\{E_{i,k}, F_{i,k}\} = O\left(\frac{\mathbf{P}\{E_{i,k}\} \cdot b(n,k)^{9/57}}{n-k}\right).$$

By the independence of disjoint sections of the random walk, $A_{i,k}$ and $B_{i,k}$ are independent of $E_{i,k}$ and $F_{i,k}$, so $\mathbf{P}\{A_{i,k}, D_{i,k}, E_{i,k}, F_{i,k}\} = \mathbf{P}\{A_{i,k}, D_{i,k}\} \cdot \mathbf{P}\{E_{i,k}, F_{i,k}\}$, and by (55) and (56) we thus have

$$(57) \qquad \mathbf{P}\{A_{i,k}, D_{i,k}, E_{i,k}, F_{i,k}\} = O\left(\frac{\mathbf{P}\{A_{i,k}\}\mathbf{P}\{E_{i,k}\} \cdot b(n,k)^{18/57}}{k(n-k)}\right).$$

Finally, the events $A_{i,k}$, $E_{i,k}$ and $\{m-1 \leq S_n \leq m\}$ all simply restrict the value of a certain random walk at its endpoint; by applying Theorem 1 to bound $\mathbf{P}\{A_{i,k}\}$, $\mathbf{P}\{E_{i,k}\}$, and $\mathbf{P}\{m-1 \leq S_n \leq m\}$, it is easily seen that

$$\mathbf{P}\{A_{i,k}\} \cdot \mathbf{P}\{E_{i,k}\} = \Theta\left(\mathbf{P}\{m-1 \leq S_n \leq m\} \cdot \frac{n^{1/2}}{(k(n-k))^{1/2}}\right).$$

Combining this last equation with (57) proves (52) and completes the proof. $\square$

It is interesting to compare Lemmas 11 and 12 with the work of Bramson [8] in studying branching Brownian motion; he required quite similar bounds on the behavior of a Brownian bridge (or equivalently, of Brownian motion conditioned on its value after some time $t$). It would be interesting to find a proof of Theorem 3, or at least of Lemmas 11 and 12, that proceeded via comparison with Brownian motion. However, such an approach may be quite difficult, given the lack of success to date at transferring Bramson's argument to a discrete setting. At any rate, with these lemmas under our belt we are ready to prove Lemma 9.

4.3. *Proof of Lemma 9.* The chain of reasoning is similar to that in the proof of Lemma 8. We first recall that for a node $v$ at depth $n$, by Lemma 7

$$(58) \qquad \mathbf{P}\{v \in G_{n,m'}\} = \Theta\left(\frac{\mathbf{P}\{m'-1 \leq S_v \leq m'\}}{n}\right).$$

Recalling the event $B_a$ defined just before Lemma 12, we will write $B_a(v)$ for the event that $B_a$ occurs for the random walk ending at $v$. Since $|\Lambda'(t)n - m'| = O(\log n)$, by Lemmas 11 and 12, for all integers $a < -1$,

$$\mathbf{P}\{m'-1 \leq S_v \leq m, S_v \text{ abo } a, B_a(v)\} = \Theta\left(\frac{\mathbf{P}\{m'-1 \leq S_v \leq m\}}{n|a|^9}\right)$$

$$= O\left(\frac{\mathbf{P}\{v \in G_{n,m'}\}}{|a|^9}\right).$$

It follows that there is some large constant $C > 0$ such that if $v \in G_{n,m'}$ then, letting $R_0, R_1, \ldots, R_n$ be the partial sums on the path to $v$, with probability at least $1/2$

$$(59) \qquad R_k \geq R_n \cdot \frac{k}{n} + b(n,k)^{1/57} \qquad \text{for all } C \leq k \leq n-C.$$

If $v \in G_{n,m'}$ and $v$ additionally satisfies (59), we say that $v$ is *well-behaved*. We write $W_{n,m'}$ for the set of well-behaved nodes in $G_{n,m'}$. We emphasize that $W_{n,m'} \subseteq G_{n,m'} \subseteq N_{n,m'}$. Furthermore, since each node in $G_{n,m'}$ is well-behaved with probability at least $1/2$, it follows from the definition of $m'$ that $\mathbf{E}|W_{n,m'}| = \Omega(\mathbf{E}|G_{n,m'}|) = \Omega(1)$. We claim that

$$\mathbf{E}\{|W_{n,m'}| \,|\, v_n \in W_{n,m'}\} = \mathbf{E}|W_{n,m'}| + O(1). \tag{60}$$

Applying (25) with $\mathcal{R} = W_{n,m'}$, it immediately follows that $\mathbf{P}\{|W_{n,m'}| \geq 1\} = \Omega(1)$, so since $W_{n,m'} \subseteq N_{n,m'}$ we have $\mathbf{P}\{|N_{n,m'}| \geq 1\} = \Omega(1)$ as claimed. We thus turn to proving (60).

We now remind the reader of some notation from earlier in the section, and introduce a few new terms. We recall that node $v_{i-1}$ has child $v_i^{(0)} = v_i$ that is on the distinguished path $P$, and that its remaining children are $v_i^{(1)}, \ldots, v_i^{(d-1)}$. Node $v_i^{(j)}$ is the root of a subtree of $T_n$ that we shall call $T_n^{i,j}$. We let $W_{n,m}^{i,j}$ be subset of $W_{n,m}^{i,j}$ contained in $T_n^{i,j}$; this is consistent with the notation $G_{n,m}^{i,j}$ and $N_{n,m}^{i,j}$ introduced in and just after (26). For each $0 \leq i \leq n-1$ fix an arbitrary node $x_i$ at depth $n$ that is a descendent of $v_{i+1}^{(1)}$, and let the partial sums of the labels on the path from $v_i$ to $x_i$ be $S_{x_0,0} = 0, S_{x_i,1}, S_{x_i,2}, \ldots, S_{x_i,n-i}$. We remark that the edge labels contributing to the sum $S_{x_i,n-i}$ are a subset of the edge labels contributing to the vertex label $S_{x_i}$; more precisely,

$$S_{x_i} = S_{x_i,n-i} + S_{v_i}. \tag{61}$$

Finally, to simplify notation, for $0 \leq i \leq n-1$ let $W_i = W_{n,m'}^{i,1}$, and denote the functions $\mathbf{P}\{\cdot | v_n \in W_{n,m'}\}$ and $\mathbf{E}\{\cdot | v_n \in W_{n,m'}\}$ by $\mathbf{P}^w\{\cdot\}$ and $\mathbf{E}^w\{\cdot\}$, respectively. We now mimic the portion of the proof of Lemma 8 that leads to (33), in our case for the particular value $m = m'$. By symmetry, we have

$$
\begin{aligned}
\mathbf{E}^w\{|W_{n,m'}|\} &= 1 + \sum_{i=1}^{n} \sum_{j=1}^{d-1} \mathbf{E}^w\{W_{n,m'}^{i,j}\} \\
&= 1 + \sum_{i=1}^{n} (d-1)\mathbf{E}^w\{W_i\}.
\end{aligned}
\tag{62}
$$

Since $d$ does not depend on $n$ and, for a given $i$, $W_{h'-i} \leq i^d$, it follows that for any fixed integer $0 < c = O(1)$, $\sum_{i=n-c}^{n}(d-1)\mathbf{E}^w\{W_i\} \leq \sum_{i=1}^{c}(d-1)i^d = O(1)$. By this fact, by (62), and since $(d-1)$ is constant, to prove the lemma it therefore suffices to show that

$$\sum_{i=1}^{h'-C} \mathbf{E}^w\{|W_i|\} = O(1), \tag{63}$$

where $C$ is the same constant as in (59).

By symmetry, $\mathbf{E}^w\{|W_i|\} = d^{n-i}\mathbf{P}^w\{x_i \in W_i\}$. In order for $x_i \in W_i$ to occur, we must in particular have that $S_{x_i} \leq m'$, so by (61), we must have $S_{x_i,n-i} \leq m' - S_{v_i}$. The fact that $v_n$ is well-behaved allows us to bound $m' - S_{v_i}$ from above. Furthermore, $S_{x_i,n-i}$ is distributed as $S_{v_{n-i}}$ and is independent of $S_{v_i}$. This independence will allow use use Theorem 1 to bound the conditional probability that $S_{x_i,n-i} \leq m' - S_{v_i}$. When $i$ is far from 1 and from $n - C$, the bounds on $S_i$ given by the fact that $v_n$ is well-behaved will ensure that $S_{v_i}$ is large enough that the conditional probability that $S_{x_i,n-i} \leq m' - S_{v_i}$ is extremely small. By slightly modifying this same approach, we will prove similar bounds when $i$ is near 1 or near $n-C$; summing these bounds will prove (63). We now turn to the details.

Fix some function $g$ with $g(k) \leq \log k$ for all $k$ and with $g(k)$ tending to plus infinity as $k$ tends to infinity. For any $i$, for all $c$ with $|c| \leq \sqrt{n-i}/g(i)$, by Theorem 1 and by the independence of $X_{v_1}, \ldots, X_{v_n}$ from $S_{x_i,n-i}$, we have

$$
\begin{aligned}
&\mathbf{P}^w\{S_{x_i,n-i} \leq \Lambda'(t)(n-i) + c\} \\
(64) \quad &= \mathbf{P}\{S_{x_i,n-i} \leq \Lambda'(t)(n-i) + c | S_{x_i,n-i} < \infty\} \cdot \mathbf{P}\{S_{x_i,n-i} < \infty\} \\
&= \Theta\left(\frac{e^{-tc}}{\sqrt{n}} \cdot \frac{1}{d^{n-i}}\right).
\end{aligned}
$$

Furthermore, by (61),

$$
(65) \quad \mathbf{P}^w\{S_{x_i} \leq m'\} = \mathbf{P}^w\{S_{x_i,n-i} \leq m' - S_i\}.
$$

Let $r = r(n) = \lceil (2\log n/|t|)^{57} \rceil$.

CASE 1 ($r \leq i \leq n-r$). For any $k \geq C$, since $v_n \in W_{n,m'}$, we have $S_{v_k} \geq m' \cdot (k/n) - 1 + b(n,k)^{1/57}$, so

$$
(66) \quad m' - S_{v_i} \leq m' \cdot \left(\frac{n-i}{n}\right) + 1 - b(n,i)^{1/57}
$$

$$
(67) \quad \leq \Lambda'(t)(n-i) + \frac{3\log n}{2|t|} - b(n,i)^{1/57}.
$$

[We will also use (66) in the case that $i < n-r$.] When $b(n,i)^{1/57} \geq 2\log n/|t|$, certainly $\log n/2|t| \leq \sqrt{n-i}/g(n-i)$, so by (64), (65) and (67) we have

$$
\mathbf{P}^w\{S_{x_i} \leq m'\} \leq \mathbf{P}^w\left\{S_{x_i,n-i} \leq \Lambda'(t)(n-i) - \frac{\log n}{2|t|} + 1\right\} = \Theta\left(\frac{1}{d^{n-i}n}\right),
$$

(68)

so by linearity of expectation and symmetry, $\mathbf{E}^w|W_i| = O(1/n)$ for such $i$. Letting $r = r(n) = \lceil (2\log n/|t|)^{57} \rceil$, then, we have

$$
\sum_{i=r}^{n-r} \mathbf{E}^w|W_i| = O\left(\frac{n-2r}{n}\right) = O(1).
$$

This bounds the bulk of the sum (63); it remains to consider the cases when $i$ is either close to 1 or close to $n - C$.

CASE 2 ($n - r \leq i \leq n - C$). Let $k = n - i$, so $C \leq k < r$. Since $k = O((\log n)^{57})$, $m' \cdot (k/n) = \Lambda'(t)k + O(1)$, so from (65), (66) and the fact that $b(n, i) = k$ we have

$$\mathbf{P}^w\{x_i \in W_{n,m'}\} \leq \mathbf{P}^w\{S_{x_i} \leq m'\} = \mathbf{P}\{S_{x_i,n-i} \leq m' - S_{v_i}\}$$
$$= O(\mathbf{P}\{S_{x_i,n-i} \leq \Lambda'(t)k - b(n,i)^{1/57}\})$$
$$= O(\mathbf{P}\{S_{x_i,n-i} \leq \Lambda'(t)k - k^{1/57}\}),$$

which together with (64) gives

$$\mathbf{P}^w\{x_i \in W_{n,m'}\} = O\left(\frac{e^{tk^{1/57}}}{d^k}\right).$$

By linearity of expectation and by symmetry, it follows that $\mathbf{E}^w|W_i| = O(e^{tk^{1/57}}) = O(k^{-3})$ (recall that $t$ is negative). Therefore

$$\sum_{i=n-r}^{n-C} \mathbf{E}^w|W_i| = O\left(\sum_{k=C}^{r} \frac{1}{k^3}\right) = O(1).$$

CASE 3 ($1 \leq i \leq r$). Let $\Delta_i = S_{v_i} - m' \cdot (i/n)$. Since $v_n \in G_{n,m'}$, necessarily $\Delta_i \geq -i/n \geq -1$. In order that $x_i \in G_{n,m'}$, it is necessary that:

- $m' \cdot ((n-i)/n) - (\Delta_i + 1) \leq S_{x_i,n-i} \leq m' \cdot ((n-i)/n) - \Delta_i$ (call this event $E_i$), and that
- $S_{x_i,n-i}$ $abo$ $-(\Delta_i + 1)$,

so

$$\mathbf{P}^w\{x_i \in W_{n,m'}\} \leq \mathbf{P}^w\{x_i \in G_{n,m'}\} \leq \mathbf{P}^w\{E_i, S_{x_i,n-i} \; abo \; -(\Delta_i + 1)\}.$$

We shall show that

$$(69) \qquad \mathbf{P}^w\{E_i, S_{x_i,n-i} \; abo \; -(\Delta_i + 1)\} = O\left(\frac{1}{d^{n-i}(n-i)}\right),$$

from which it follows just as above that $\sum_{i=1}^{r} \mathbf{E}^w|W_i| = O(1)$. It thus remains to prove (69). We first observe that by the same argument used to prove (68), we have

$$\mathbf{P}^w\left\{E_i, S_{x_i,n-i} \; abo \; -(\Delta_i + 1)\Big|\Delta_i \geq \frac{2\log n}{|t|}\right\}$$
$$\leq \mathbf{P}^w\left\{E_i\Big|\Delta_i \geq \frac{2\log n}{|t|}\right\}$$

(70)
$$\leq \mathbf{P}\left\{ S_{x_i, n-i} \leq \Lambda'(t)(n-i) - \frac{\log n}{2|t|} \right\}$$

$$= O\left( \frac{1}{d^{n-i} \cdot n} \right).$$

Furthermore, if $\Delta_i < 2\log n/|t|$, then $\Delta_i \leq \sqrt{n-i}/g(n-i)$, so by applying first Lemma 11, then Theorem 1, we obtain

$$\mathbf{P}^w\left\{ E_i, S_{x_i, n-i} \ abo \ -(\Delta_i + 1) \Big| \Delta_i < \frac{2\log n}{|t|} \right\}$$

(71)
$$= O\left( \frac{\Delta_i^9 \mathbf{P}^w\{E_i | \Delta_i < (2\log n)/|t|\}}{n-i} \right)$$

$$= O\left( \frac{\Delta_i^9 e^{\Delta_i t}}{n-i} \right)$$

$$= O\left( \frac{1}{d^{n-i}(n-i)} \right).$$

Combining (70) and (72) proves (69) and completes the proof of Lemma 9.

**5. The lower bound.** We know that the naive approach to proving a lower bound on $\mathbf{P}\{M_n \leq m\}$, namely, bounding $\mathbf{P}\{v_n \in N_{n,m}\}$, then applying a union bound, will not work. It is the approach we used in Section 1.1 when discussing many *independent* random walks, and only begins to yield tail bounds when $m \leq m^* + O(1)$.

On the other hand, we observe that we can *easily* prove strong enough bounds for one group of potential nodes of $N_{n,m}$ when $m \leq m' + O(1)$, namely, the set $G_{n,m}$ [we recall that $m' = \Lambda'(t)n - (3\log n)/2t$]. By Lemma 11 (or, indeed, by Lemma 7), when $m$ is not too far from $m'$, $\mathbf{P}\{v_h \in G_{n,m}\} = O(\mathbf{P}\{v_h \in N_{n,m}\}/n)$, from which strong bounds on $\mathbf{P}\{G_{n,m} \neq \varnothing\}$ follow directly from Theorem 1 and a union bound. Intuitively, this bound on $G_{n,m}$ contains almost the whole "reason for" the lower bound, in the following sense. Given a node $v$ at depth $n$ and $m \leq m'$, if $S_v \leq m$ is to occur, then *either* $S_v$ is "close to being" a leading node, in the sense that there is a "small" constant $a$ for which $S_v$ stays above $-a$, or $S_v$ is "far from being" a leading node, in which case there is a "large" constant $c$ for which $S_v$ does *not* stay above $-c$ (or something in between these two scenarios occurs). The former scenario is unlikely due to the bounds in Lemma 11, and the latter scenario is unlikely due to large deviations estimates (and everything in between is unlikely due to a mix of these two reasons).

We prove our general bound by splitting $N_{n,m}$ into many groups which differentiate among the possibilities outlined in the above sketch. We then bound the probability each group is nonempty using large deviations estimates and the Lemmas from Section 4.2. When we recombine the bounds we obtain for the individual groups (by union bound), the result will be strong enough to yield bounds on $\mathbf{P}\{N_{n,m} \neq \varnothing\}$ that are exponentially small in $m' - m$ when $m$ is not too far from $m'$. For convenience, we restate the lower bound we are aiming to prove.

LEMMA 13. *There are $C_2 > 0$, $\delta_2 > 0$ such that for all $x > 0$ we have* $\mathbf{P}\{M_n \leq m' - x\} \leq C_2 e^{-\delta_2 x}$.

Theorem 3 follows easily from the above lemma and Corollary 10; its proof appears in Section 6, below. In proving Lemma 13, it will be useful to use Chernoff [14] bounds as well as Theorem 1, as Chernoff bounds are not limited to windows of size $o(\sqrt{n})$ around a linear slope. In particular, for our choice of $t$ we have:

LEMMA 14 (Chernoff bound). *For all $r > 0$ and all integers $k \geq 1$,*

$$\mathbf{P}\{S_{v_k} \leq \Lambda'(t)k - r\} \leq \frac{e^{tr}}{d^k}.$$

Lemma 13 is a rather straightforward consequence of the following, weaker lemma, plus Lemma 14:

LEMMA 15. *There are $C_3 > 0$, $\delta_3 > 0$ such that for all $m < m'$ with* $m' - m \leq (2\log n)/|t|$,

$$\text{(72)} \qquad \mathbf{P}\{\exists v \in N_n : m - 1 \leq S_v \leq m\} \leq C_3 e^{-\delta_3(m'-m)}.$$

PROOF OF LEMMA 13 ASSUMING LEMMA 15. Given $m$ with $0 \leq m' - m \leq \sqrt{n/g(n)}$, if $m' - m \geq 2\log n/|t|$ then, letting $\Delta = \Lambda'(t)n - m$, by the definition of $m'$ we have $\Delta \geq (m' - m)/4$. It follows by a union bound and by Theorem 1 that

$$\mathbf{P}\{M_n \leq m\} \leq d^n \mathbf{P}\{S_{v_n} \leq \Lambda'(t)n - \Delta\} = O(e^{t\Delta}) = O(e^{t(m'-m)/4}),$$

which proves Lemma 13 in the case that $m' - m \geq (2\log n)/|t|$ (since $t$ is negative). If $0 \leq m' - m < (2\log n)/|t|$, then we have the following inclusion:

$$
\text{(73)} \qquad
\begin{aligned}
\{M_n \leq m\} &\subseteq \left\{ M_n \leq m' - \frac{2\log n}{|t|} \right\} \\
&\quad \cup \bigcup_{i=\lfloor m'-m \rfloor}^{\lfloor (2\log n)/|t| \rfloor} \{\exists v \in N_n : m' - i - 1 \leq S_v \leq m' - i\}.
\end{aligned}
$$

Since $m' - (2 \log n)/|t| = \Lambda'(t)n - (\log n)/2|t|$, by Lemma 14 and a union bound,

$$(74) \qquad \mathbf{P}\left\{M_n \le m' - \frac{2 \log n}{|t|}\right\} \le \frac{1}{\sqrt{n}} = e^{-\log n/2} \le e^{t(m'-m)/4}.$$

By (73), (74), Lemma 15 and a union bound, we thus have

$$\mathbf{P}\{M_n \le m\} \le e^{t(m'-m)/4}$$
$$+ \sum_{i=\lfloor m'-m \rfloor}^{\lfloor (2 \log n)/|t| \rfloor} \mathbf{P}\{\exists v \in N_n : m' - i - 1 \le S_v \le m' - i\}$$
$$\le e^{t(m'-m)/4} + C_3 \sum_{i=\lfloor m'-m \rfloor}^{\lfloor (2 \log n)/|t| \rfloor} e^{-\delta_3 i}$$
$$\le e^{t(m'-m)/4} + C_3 e^{-\delta_3 \lfloor m'-m \rfloor} \left(\frac{1}{1 - e^{-\delta_3}}\right),$$

which proves Lemma 13 in the case that $0 \le m' - m \le (2 \log n)/|t|$.  $\square$

We now turn our attention to proving Lemma 15. We first note that we need only prove (72) for $m' - m$ larger than any fixed constant, as we may presume the bound holds for $m' - m = O(1)$ by our choice of $C_3$. For the remainder of the section, we write $\Delta = \Lambda'(t)n - m$ and $\Delta' = m' - m > L$ for some large constant $L$. We have $\Delta' = \Delta + (3 \log n)/2|t|$, so $\Delta' \le (2 \log n)/|t|$ and $\Delta \le (\log n)/2|t|$.

We now proceed to define the "homogeneous" groups discussed above, and bound the probabilities they are nonempty, in a sequence of claims. The proof of each claim will consist of straightforward applications of Theorem 1 and Lemmas 11 and 12, and Lemma 15 will be a trivial consequence of the bounds of the claims.

Let $a = \lfloor e^{|t| \Delta'/10} \rfloor$, and let $A_a$ be the set of nodes $v \in N_n$ for which $m - 1 \le S_v \le m$ and for which $\{S_v \ abo \ - a\}$ occurs.

CLAIM 16.    $\mathbf{P}\{A_a \ne \varnothing\} = O(e^{9t(m'-m)/10})$.

PROOF.    By Lemma 11,

$$\mathbf{P}\{v_n \in A_a | m - 1 \le S_v \le m\} = O\left(\frac{a^9}{n}\right) = O\left(\frac{e^{9t(m'-m)/10}}{n}\right).$$

By Theorem 1, it follows that

$$\mathbf{P}\{v_n \in A_a\} = O\left(\frac{e^{9t(m'-m)/10}}{d^n}\right).$$

The claim follows by a union bound over all $v \in N_n$. $\square$

Next, let $c = \lfloor \log^{10} n \rfloor$, and let $H_c$ consist of the nodes $v \in N_n$ with $m-1 \leq S_v \leq m$ and for which $\{S_v \ abo \ -c\}$ *does not* occur.

CLAIM 17. $\mathbf{P}\{H_c \neq \varnothing\} = O(e^{-(m'-m)})$.

PROOF. If $H_c \neq \varnothing$ then for some $0 < k < n$ there is a node $y \in N_k$ for which $S_y \leq mk/n - c$. We have

$$\frac{mk}{n} - c = \frac{(\Lambda'(t)n - \Delta)k}{n} - c = \Lambda'(t)k - \frac{\Delta k}{n} - c.$$

Therefore, by Lemma 14,

$$\tag{75} \mathbf{P}\left\{S_{v_k} \leq \frac{mk}{n} - c\right\} = O\left(\frac{e^{t(c+\Delta k/n)}}{d^k}\right).$$

Since $c = \lfloor \log^{10} n \rfloor$, we may assume $n$ is large enough that $c - \Delta k/n \geq \log^9 n$, so (75) yields

$$\tag{76} \mathbf{P}\left\{S_{v_k} \leq \frac{mk}{n} - c\right\} = O\left(\frac{e^{t \log^9 n}}{d^k}\right).$$

By a union bound, we obtain

$$\mathbf{P}\{H_c \neq \varnothing\} \leq \sum_{k=1}^{n-1} d^k \mathbf{P}\left\{S_{v_k} \leq \frac{mk}{n} - c\right\}$$

$$= O\left(\sum_{k=1}^{n-1} e^{t \log^9 n}\right) = O(e^{\log n + t \log^9 n}) = O(e^{-(m'-m)}). \quad \square$$

For each integer $b$ with $a \leq b \leq c$, we let $M_b$ be the set of vertices $v \in N_n$ for which $m - 1 \leq S_v \leq m$ and for which $\{S_v \ abo \ -(b+1)\}$ *occurs* but $\{S_v \ abo \ -(b+1)\}$ *does not occur*. We next define a subset of $M_b$ which we call $M_b^{\mathrm{mid}}$ by saying that $v_n \in M_b^{\mathrm{mid}}$ if $v_n \in M_b$ and additionally, there is $k$ with $b^{57} \leq k \leq n - b^{57}$ for which $S_{v_k} \leq mk/n - b$. We extend this definition to all nodes $v \in N_n$ by symmetry.

Similarly, for values $k$ with $k < b^{57}$ or with $n - k < b^{57}$, we define a set $M_b^k$ by saying that $v_n \in M_b^k$ if $v_n \in M_b$ and additionally, $k$ is the *smallest value* for which $S_{v_k} \leq mk/n - b$; again, we extend this definition to all nodes $v \in N_n$ by symmetry. The sets $M_b^{\mathrm{mid}}$ and $\{M_b^k : \min(k, n-k) < b^{57}\}$ partition $M_b$. We now bound the probabilities that these sets are nonempty.

CLAIM 18. $\mathbf{P}\{\bigcup_{b=a}^c \{M_b^{\mathrm{mid}} \neq \varnothing\}\} = O(e^{t(m'-m)/11})$.

Proof.   Fix an integer $b$ with $a \le b \le c$. By Lemma 12,

$$\mathbf{P}\{v_n \in M_b^{\mathrm{mid}} | m-1 \le S_{v_n} \le m\} = O\left(\frac{1}{nb^9}\right),$$

so by Theorem 1,

$$\mathbf{P}\{v_n \in M_b^{\mathrm{mid}}\} = O\left(\frac{\mathbf{P}\{m-1 \le S_{v_n} \le m\}}{nb^9}\right)$$

$$= O\left(\frac{ne^{t(m'-m)}}{nb^9 d^n}\right)$$

$$= O\left(\frac{e^{t(m'-m)}}{b^9 d^n}\right).$$

By a union bound, it follows that

$$\mathbf{P}\{M_b^{\mathrm{mid}} \ne \varnothing\} = O\left(\frac{e^{t(m'-m)}}{b^9}\right),$$

so by summing over $b$ with $a \le b \le c$, since $a = e^{|t|(m'-m)/11}$ we obtain

$$\mathbf{P}\left\{\bigcup_{b=a}^c \{M_b^{\mathrm{mid}} \ne \varnothing\}\right\} = O\left(\sum_{b=a}^c \frac{e^{t(m'-m)}}{b^9}\right) = O\left(\frac{e^{t(m'-m)}}{a^{10}}\right) = O(e^{t(m'-m)/11}). \quad \square$$

Claim 19.   $\mathbf{P}\{\bigcup_{b=a}^c \bigcup_{k=1}^{b^{57}-1} \{M_b^k \ne \varnothing\}\} = O(e^{-(m'-m)})$.

Proof.   Fix $b$ and $k$ as above. For each node $x$ at depth $k$, let $W_x$ be the set of descendents of $x$ in $M_b^k$. By a union bound over $N_k$, $\mathbf{P}\{M_b^k \ne \varnothing\} \le d^k \mathbf{P}\{W_{v_k} \ne \varnothing\}$. If $W_{v_k}$ is nonempty, then necessarily $mk/n - (b+1) \le S_{v_k} \le mk/n - b$. Since $k \le b^{57} \le c^{57} \le (\log n)^{171}$,

$$\frac{mk}{n} - b = \Lambda'(t)k + \frac{\Delta k}{n} - b = \Lambda'(t)k - b + o(1).$$

It follows by Lemma 14 that

$$(77) \qquad \mathbf{P}\{W_{v_k} \ne \varnothing\} = \Theta(\mathbf{P}\{S_{v_k} \le \Lambda'(t)k - b\})$$

$$= O\left(\frac{e^{tb}}{d^k}\right).$$

By (77) and a union bound, we have that $\mathbf{P}\{M_b^k \ne \varnothing\} = O(e^{tb})$, so by a second union bound

$$\mathbf{P}\left\{\bigcup_{b=a}^c \bigcup_{k=1}^{b^{57}-1} \{M_b^k \ne \varnothing\}\right\} = O\left(\sum_{b=a}^c b^{57} e^{tb}\right) = O(a^{57} e^{ta}) = O(e^{-(m'-m)}). \quad \square$$

CLAIM 20. $\mathbf{P}\{\bigcup_{b=a}^{c}\bigcup_{k=1}^{b^{57}-1}\{M_b^{n-k}\neq\varnothing\}\}=O(e^{-(m'-m)})$.

PROOF. Fix $b$ and $k$ as above. For each node $x$ at depth $n-k$, let $W_x$ be the set of descendants of $x$ in $M_b^{n-k}$. By a union bound,

$$(78) \qquad \mathbf{P}\{M_b^{n-k}\neq\varnothing\}\leq d^{n-k}\mathbf{P}\{W_{v_{n-k}}\neq\varnothing\}.$$

Suppose $W_{v_r}$ is nonempty—then necessarily

$$\frac{m(n-k)}{n}-(b+1)\leq S_{v_{n-k}}\leq\frac{m(n-k)}{n}-b,$$

*and in addition* $\{S_{v_{n-k}}\ abo\ -b\}$ must occur (or else $n-k$ is not the *first* time the random walk ending at $v_n$ falls $b$ below its conditional mean). Now, by Theorem 1, since $\Delta\geq(3\log n)/2t$,

$$(79) \qquad \begin{aligned}
\mathbf{P}\left\{S_{v_{n-k}}\leq\frac{m(n-k)}{n}-b\right\}&=\mathbf{P}\left\{S_{v_{n-k}}\leq\Lambda'(t)(n-k)-\frac{\Delta(n-k)}{n}-b\right\}\\
&=\Theta\left(\frac{e^{t(\Delta(n-k)/n+b)}}{d^{n-k}\sqrt{n-k}}\right)=\Theta\left(\frac{e^{t(\Delta+b)}}{d^{n-k}\sqrt{n}}\right)\\
&=O\left(\frac{ne^{tb}}{d^{n-k}}\right).
\end{aligned}$$

Furthermore, by Lemma 11,

$$(80) \qquad \begin{aligned}
\mathbf{P}\left\{S_{v_{n-k}}\ abo\ -b\,\Big|\,\frac{m(n-k)}{n}-(b+1)\leq S_{v_{n-k}}\leq\frac{m(n-k)}{n}-b\right\}\\
=O\left(\frac{b^9}{n-k}\right)=O\left(\frac{b^9}{n}\right).
\end{aligned}$$

Combining (79) and (80) yields that

$$\mathbf{P}\{W_{v_{n-k}}\neq\varnothing\}=O\left(\frac{b^9e^{tb}}{d^{n-k}}\right)=O\left(\frac{e^{9\log b+tb}}{d^{n-k}}\right),$$

which combined with (78) implies that $\mathbf{P}\{M_b^{n-k}\neq\varnothing\}=O(e^{9\log b+tb})$. Just as in the proof of Claim 19, summing this bound over $b$ and $k$ yields the result. $\square$

We are now prepared for:

PROOF OF LEMMA 15. It is immediate from the definitions of the sets $A_a$, $H_c$, $M_b^{\mathrm{mid}}$ and $M_b^k$ that any vertex $v\in N_n$ with $m-1\leq S_v\leq m$ is either in $A_a$, or in $H_c$, or in $M_b^{\mathrm{mid}}$ (for some integer $a\leq b\leq c$), or in $M_b^k$ (for some integer $a\leq b\leq c$ and some integer $k$ for which either $1\leq k\leq b^{33}$ or $h - b^{33} \leq k \leq h$). Applying Claims 16–20, respectively, to bound each of these events, it follows that

$$\mathbf{P}\{\exists v \in N_n : m - 1 \leq S_v \leq m\} = O(e^{9t(m'-m)/10}) + O(e^{-(m'-m)})$$
$$+ O(e^{t(m'-m)/11}) + O(e^{-(m'-m)})$$
$$+ O(e^{-(m'-m)})$$
$$= O(e^{-\min(|t|/11,1)(m'-m)}),$$

so the bound of Lemma 15 holds as long as we choose $\delta_1$ so that $0 < \delta_1 \leq \min(|t|/11, 1)$ and choose $C_1$ large enough.   $\square$

**6. Proof of Theorem 3.**  By Corollary 10, there exist $C_1 > 0$, $\delta_1 > 0$ such that for all $x > 0$,

$$\mathbf{P}\{M_n \geq m' + x | \mathcal{S}\} \leq C_1 e^{-\delta_1 x}.$$

By Lemma 13 and Bayes' formula, there exist $C_2 > 0$, $\delta_2 > 0$ such that for all $x > 0$,

$$\mathbf{P}\{M_n \leq m' - x | \mathcal{S}\} \leq \frac{\mathbf{P}\{M_n \leq m' - x\}}{\mathbf{P}\{\mathcal{S}\}} \leq \frac{C_2}{\mathbf{P}\{\mathcal{S}\}} e^{-\delta_2 x}.$$

Taking $C_4 = \max\{C_1, C_2/\mathbf{P}\{\mathcal{S}\}\}$ and $\delta = \min\{\delta_1, \delta_2\}$, we obtain that for all $x > 0$,

$$\mathbf{P}\{|M_n - m'| \geq x | \mathcal{S}\} \leq C_4 e^{-\delta x}.$$

It follows immediately that

$$|\mathbf{E}\{M_n | \mathcal{S}\} - m'| \leq \mathbf{E}\{|M_n - m'| \,|\mathcal{S}\} \leq \sum_{i=0}^{\infty} \mathbf{P}\{|M_n - m'| \geq i | \mathcal{S}\}$$
$$\leq \frac{C_4}{1 - e^{-\delta}},$$

proving (7); and letting $\Delta = C_4/(1 - e^{-\delta})$, for all $x > 0$ we have

$$\mathbf{P}\{|M_n - \mathbf{E}\{M_n | \mathcal{S}\}| > x | \mathcal{S}\} \leq \mathbf{P}\{|M_n - m'| > x - \Delta | \mathcal{S}\} \leq C_4 e^{\Delta \delta} e^{-\delta x},$$

so (8) holds with this choice of $\delta$ and with $C = C_4 e^{\Delta \delta}$.

**Acknowledgment.**  The authors would like to thank an anonymous referee for many useful comments and suggestions.

DÉPARTEMENT DE MATHEMATIQUES
ET DE STATISTIQUE
UNIVERSITÉ DE MONTRÉAL
C.P. 6138, SUCC. CENTRE-VILLE
MONTRÉAL, QC, H3C 3J7
CANADA
E-MAIL: addario@dms.umontreal.ca

CANADA RESEARCH CHAIR
SCHOOL OF COMPUTER SCIENCE
MCGILL UNIVERSITY
3480 UNIVERSITY STREET
MONTRÉAL, QC, H3A 2A7
CANADA
AND
EQUIPE MASCOTTE, LABO. I3S
INRIA
SOPHIA ANTIPOLIS
FRANCE
E-MAIL: breed@cs.mcgill.ca